\documentclass[12pt]{article}
\usepackage{amsmath,amssymb} 
\usepackage{graphicx} 
\usepackage{color}
\usepackage{caption}
\usepackage{subcaption}
\newtheorem{theorem}{Theorem}[section]
\newtheorem{lemma}[theorem]{Lemma}

\newtheorem{proposition}[theorem]{Proposition}
\newtheorem{definition}[theorem]{Definition}
\newtheorem{assumption}[theorem]{Assumption}
\newtheorem{remark}{Remark}
\newtheorem{example}{Example}
\def\ulI{\underline{I}}
\def\calA{\mathcal{A}}
\def\calG{\mathcal{G}}
\def\eqref#1{(\ref{#1})}
\def\tfrac#1#2{{\textstyle\frac{#1}{#2}}}
\def\mod{\mathbb{M}}
\def\obs{\mathbb{O}}
\def\tint{{\textstyle\int}}
\begin{document}

\title{All-at-once formulation meets the Bayesian approach: A study of two prototypical linear inverse problems}
\author{Anna Schlintl and Barbara Kaltenbacher\\
Alpen-Adria-Universit\"at Klagenfurt, 
Austria\\
barbara.kaltenbacher@aau.at, anna.schlintl@aau.at
}
\maketitle

\begin{abstract}In this work, the Bayesian approach to inverse problems is formulated in an all-at-once setting. 
The advantages of the all-at-once formulation are known to include the avoidance of a parameter-to-state map as well as numerical improvements, especially when considering nonlinear problems. In the Bayesian approach, prior knowledge is taken into account with the help of a prior distribution. In addition, the error in the observation equation is formulated by means of a distribution. This method naturally results in a whole posterior distribution for the unknown target, not just point estimates. This allows for further statistical analysis including the computation of credible intervals. 
We combine the Bayesian setting with the all-at-once formulation, resulting in a novel approach for investigating inverse problems. With this combination we are able to chose a prior not only for the parameter, but also for the state variable, which directly influences the parameter. Furthermore, errors not only in the observation equation,
but additionally, in the model can be taken into account. 
We analyze this approach with the help of two linear standard examples, namely the inverse source problem for the Poisson equation and the backwards heat equation, i.e. a stationary and a time dependent problem. Appropriate function spaces and derivation of adjoint operators are investigated. To assess the degree of ill-posedness, we analyze the singular values of the corresponding all-at-once forward operators.
Finally, joint priors are designed and numerically tested. 
\end{abstract}

%
\vspace{2pc}
\noindent{\it Keywords}: inverse problems, all-at-once formulations, Bayesian inverse problems
\section{Introduction.}
\normalsize
In a general setting, the model to work with in the Bayesian approach is given by
\begin{equation}
\label{modeleq}
y^{\delta} = G x + \delta \eta,
\end{equation}
where $G$ is the (linear) forward operator, mapping between Hilbert spaces $X$ and $Y$, $\delta\geq0$ describes the noise level, $x$ is the unknown target and $y^{\delta}$ denotes the data, when $\eta$ is a random variable describing the noise. We restrict our considerations to Gaussian noise, i.e. $\eta \sim \mathcal{N}(0, \Sigma)$. This equation can be transformed such that it results in a problem under Gaussian white noise given by
\begin{equation}
z^{\delta} = \Sigma^{-1}G x + \delta \Sigma^{-1} \eta.
\end{equation}
The prior distribution for $x$ is chosen to be normal, i.e. $\eta \sim \mathcal{N}(0, \frac{\delta^2}{\alpha} C_0)$, with the noise level $\delta$ and a scaling parameter $\alpha$. In that setting, much is known, especially that the posterior is also normally distributed with mean and covariance given by
\begin{eqnarray}
x_{\alpha}^{\delta} &= C_0^{1/2}(\alpha I + H)^{-1} B^* z^{\delta}, \label{xaldel}\\
C_{\alpha}^{\delta} &= \delta^2C_0^{1/2}(\alpha I  + H)^{-1}C_0^{1/2}, \label{Caldel}
\end{eqnarray}
where $B = \Sigma^{-1/2}G C_0^{1/2}$ and  $H = B B^*$, see, e.g., (\cite{alen:bayes,  knapik:inverseheat, gugu:inverseheat, flath:inversebayes, martin:bayes, thanh:inversebayes, cotterstuartPDE, stuart:dashti, agapioumathe2018, mathe2018}.
Here and below the superscript ${}^*$ denotes the Hilbert space adjoint. 

In equation (\ref{modeleq}) the inverse problem is given in its reduced form, i.e. the model
\begin{equation}
\label{model}
\mod(u,x) = 0
\end{equation}
and the observation equation 
\begin{equation}
\label{obseq}
\obs(u) = y
\end{equation}
with operators $\mod:U\times X\to W'\times Y$, $\obs:U\to Y$ mapping between Hilbert spaces $U,X,W',Y$,
are combined through the parameter-to-state map $S$, which maps the parameter $x$ to the state $u$ and is determined by the identity
\begin{equation*}
\mod(S(x),x) = 0.
\end{equation*}
The problem then transforms to 
\[
G(x):=\obs(S(x)) = y.
\]
Contrary to the reduced approch, the all-at-once formulation combines the model (\ref{model}) and the observation equation (\ref{obseq}) in one system given by
\begin{equation}
\label{aaomodel}
\calG (\mathbf{x}) = \calG (u,x) = \left( \begin{array}{c} \mod(u,x) \\ \obs(u) \end{array} \right) = \left( \begin{array}{c} 0 \\ y \end{array} \right) = \mathbf{y}.
\end{equation}
This approach has recently gained attention, see, e.g., (\cite{burger12002, burger22002, Haber2001, kaltenbacher:aao2, kaltenbacher:aao, kaltenbacher:aao3, tramaao}).
While these paper remain in a purely deterministic setting, it is the aim of this paper to apply the Bayesian approach to the all-at once formulation.
Incorporating a priori information not only on the parameter but also on the state is expected to potentially improve reconstructions. 
Moreover, considering possible perturbations in both model and observations, i.e., assuming $\mathbf{y}^\delta = \mathbf{y}+ \delta \left( \begin{array}{c} \eta_1 \\ \eta_2 \end{array} \right)$, allows to take into account not only noise in the data but also uncertainty in the model, which is relevant in many applications. 
Combinations to a posterior distribution obtained by Bayes' formula allows to quantify the resulting uncertainty in both state and parameter.

\subsection{Examples.} \label{subs:ex}
To illustrate the all-at-once setting for inverse problems, two examples of linear inverse problems for PDEs are given in this section. These examples will serve as prototypical showcases in the analytic and numerical considerations of the following sections.
Here we will only sketch these examples; the function spaces and operator definitions will be made clear in the next section.

To avoid confusion with the space variable, we will rename the parameter $\theta$ here.  

\begin{example} \label{ex:is_intro}
The first example is a linear inverse source problem where we aim to recover the source $u$ from the model
\begin{equation}\label{invs}
\begin{array}{rcll}
- \Delta u&=& f+\theta&  \quad \mbox{ in } \Omega\\
u&=& 0& \quad \mbox{ on } \partial \Omega
\end{array}
\end{equation}
with noisy observations $y^{\delta}$ given through the equation
\begin{equation}
y^{\delta}= \obs u + \delta \eta \quad \in \Omega.
\end{equation}
The problem can be reformulated in an all-at-once fashion resulting in
\begin{eqnarray}
\label{invsaao}
\mathbf{y}^\delta = \calG \mathbf{x} + \delta \pmb{\eta} \Longleftrightarrow \left( \begin{array}{c} y^\delta_1(x) \\ y^\delta_2(x) \end{array} \right) = \left( \begin{array}{cc} - \Delta & -I \\ \obs & 0 \end{array} \right) \left( \begin{array}{c} u(x) \\ \theta(x) \end{array} \right) + \delta \left( \begin{array}{c} \eta_1 \\ \eta_2 \end{array} \right),
\end{eqnarray}
where $\mathbf{y}$ is now a vector consisting of the known inhomogeneity term $f$ and of the observation $y$, $\calG $ is a block matrix operator, the target $\mathbf{x}$ is a vector of the state variable $u$ and the source $\theta$. Additionally, the noise consists of two components, where the possible perturbation in the system is modeled through $\eta_1$ and the error in the observation equation by $\eta_2$. 

Inverse source problems with possibly more complicated elliptic operators in place of $-\Delta$ as well as restricted measurements $\obs u$, arise in numerous applications. For example, with the Helmholtz equation in place of the Poisson equation as a model, the problem corresponds to the frequency domain formulation of reconstructing sound sources from measurements taken by a microphone array, see, e.g. \cite{KKG18,Schumacher03}.
\end{example}

\begin{example} \label{ex:bh_intro}
The second example of interest in this paper is the backward heat equation. The ambition is to recover the initial temperature when measurements are only available at some point later in time. This means, contrary to the first example this problem is time dependent. We consider the model as
\begin{equation*}
\begin{array}{rcll}
\partial_t u - \Delta u &=& f \quad &\mbox{in } (0,T)\times\Omega \\
u &=& 0 \quad &\mbox{on } (0,T)\times\partial \Omega \\ 
u(0,x) &=& \theta(x) \quad &x\in \Omega.
\end{array}
\end{equation*}
with given observation equation
\begin{equation*}
y^{\delta} (x) = u(T,x) + \delta \eta \quad \mbox{ in } \Omega.
\end{equation*}
To formulate the model and the observation equation in one system, some rewriting has to be carried out first. To symbolize the observation of $v(t,x)$ in point $T$, we use the Dirac operator defined by $\delta_T v(t,x) := v(T,x)$, which maps the state to its value at the final time $T$. Further, we want to incorporate the initial condition into the state space, which we do by making the ansatz $u(t,x) := \hat{u}(t,x) + \theta(x)$ for some $\hat{u}$ contained in a linear space of functions vanishing at initial time. This (after skipping the hat) leads to the all-at-once formulation
\begin{eqnarray}
\label{isaao}
\mathbf{y}^\delta = \calG \mathbf{x} + \delta \pmb{\eta} \Longleftrightarrow\left( \begin{array}{c} y^\delta_1 \\ y^\delta_2 \end{array} \right) = \left( \begin{array}{cc} \partial_t - \Delta & -\ulI\Delta \\ \delta_T & \ulI \end{array} \right) \left( \begin{array}{c} u \\ \theta \end{array} \right) + \delta \left( \begin{array}{c} \eta_1 \\ \eta_2 \end{array} \right).
\end{eqnarray}
Here $y_1(t,x)=f(t,x)$, $y_2(x)=y(x)$ and the operator $\ulI$ maps an only space-dependent function to a formally space and time dependent function by assigning $v$ to the time-constant function $t\mapsto v$. 

A classical application that can be modeled by backwards diffusion is deconvolution of images.
Further related application examples are the identification of
airborne contaminants \cite{Akceliketal2005} and imaging with acoustic or elastic waves
in the presence of strong attenuation, arising, e.g., 
in  photoacoustic tomography \cite{KowSch12}.

\end{example}

The remainder of this paper is organized as follows.
In Section \ref{sec:adj} we provide more details on appropriate function space settings for the two prototypical examples above. In particular, since the adjoint of $\calG $ is required for the computation of the posterior mean and covariance \eqref{xaldel}, \eqref{Caldel}, but also for many other reconstruction methods, we will provide details on this. 
For both of these examples, the degree of ill-posedness, i.e., the decay reate of the eigenvalues of $G^*G$ is well known in the reduced setting and the question arises whether this behaviour may change in the all-at-once formulation. We therefore investigate the eigenvalues of the operators $\calG ^*\calG $ both analytically and numerically in Section 
\ref{sec:eig}. 
Section \ref{sec:conv} is devoted to a convergence analysis of the Bayesian all-at-once approach, where we can heavily rely on existing literature, especially the results from \cite{agapioumathe2018, mathe2018}, that largely carry over to the all-at-once formulation; however, certain conditions need a different interpretation than in the reduced setting, which we do by means of the prototypical examples above.
The important question of how to choose priors not only for the parameter but also for the state is discussed in Section \ref{sec:priors}.
Finally, in Section \ref{sec:num}, we provide some numerical results.

\section{Function space setting and computation of adjoints.} \label{sec:adj} 
For the following theoretical and numerical considerations, the two examples given in Section \ref{subs:ex} are stated more precisely in terms of functions spaces and their adjoint operators are computed in this section.
To this end, we will restrict ourselves to full observations $(\obs u)(x)=u(x)$, $x\in\Omega$, and vanishing inhomogeneity $f$, and employ the following notations.
The superscript ${}^*$ denotes the Hilbert space adjoint $B^*:Z\to A $ of a linear operator $B:A\to Z$, whereas $B^\star:Z'\to A'$ denotes its Banach space adjoint, mapping beween the dual spaces $A'$, $Z'$. 
The Laplacian $-\Delta$ with homogeneous Dirichlet boundary conditions will be denoted by $\calA $ both when acting from $H_0^1(\Omega)$ to its dual $H^{-1}(\Omega)$ and when acting from $H_0^1(\Omega) \cap H^2(\Omega)$ to $L^2(\Omega)$.
In the context of the time dependent model of the backwards diffusion problem, we will denote by $C(0,T;H)$, $L^2(0,T;H)$, $H^1(0,T;H)$ or simply shorthand $C(H)$, $L^2(H)$, $H^1(H)$, the Bochner spaces of time dependent functions with values in the space $H$ (which will typically be a space of $x$ dependent functions).
Moreover, to map elements of such a space $H$ into formally time dependent functions, we will use the operator $\ulI$ defined by 
\[ 
\ulI: H\to L^2(0,T;H)\, \quad v\mapsto (t\mapsto v)
\]
(which can as well be considered as an operator mapping into $C(0,T;H)$, or $H^1(0,T;H)$ or even $C^\infty(0,T;H)$), 
whose $H$-$L^2(0,T;H)$ Hilbert space adjoint is the averaging operator defined by
\[ 
\ulI^*: L^2(0,T;H)\to H\, \quad z\mapsto \int_0^T z(t)\, dt\,.
\]
Finally, by $(e^{-\calA t})_{t>0}$ we denote the semigroup associated with the heat equation, cf., e.g., \cite{EngelNagel99} and Subsection \ref{sec:semigroup-prior} below.

\subsection{Inverse source problem.} \label{sec:is}
The all-at-once formulation of the inverse source problem, as given in (\ref{invsaao}) consists of the vector of functions $\mathbf{x}(x) = (u(x),\theta(x))^T$ and the block operator matrix 
\begin{equation}
\calG  = \left( \begin{array}{cc} \calA & -I \\ I & 0 \end{array} \right),
\end{equation}
which acts from the space $U\times X:= H_0^1(\Omega) \cap H^2(\Omega) \times L^2(\Omega)$ to $W'\times Y:=L^2(\Omega) \times L^2(\Omega)$. 
Here $I$ denotes both the embedding of $H_0^1(\Omega)\cap H^2(\Omega)$ into $L^2(\Omega)$ and the identity on $L^2(\Omega)$. The scalar product on $U$ for elements $u,v \in U$ is defined as
\begin{equation}
\langle v(x),u(x)\rangle_U := \int_{\Omega} \calA v(x) \calA u(x) \mbox{ dx}.
\end{equation}
For functions $\mathbf{x_1}(x)=(u_1(x),\theta_1(x))$ and $\mathbf{x_2}(x)=(u_2(x),\theta_2(x)) \in U \times L^2(\Omega)$ it holds that
\begin{equation}
\langle\mathbf{x_1}, \mathbf{x_2}\rangle_{U \times L^2(\Omega)} = \int_{\Omega} \left[ \calA  u_1(x) \calA u_2(x) + \theta_1(x) \theta_2(x) \right] \mbox{ dx }.
\end{equation}
Therefore, the adjoint $\calG ^*$ of $\calG $ can be computed from
\begin{align}
\langle\calG \mathbf{x_1}, \mathbf{x_2}\rangle_{L^2 \times L^2} &= \int_{\Omega}\left[ (\calA u_1  - \theta_1) u_2 + u_1 \theta_2 \right] \mbox{ dx } \\
&= \int_{\Omega} \left[ \calA u_1 \calA (\calA ^{-1}u_2  + \calA ^{-2}\theta_2) - \theta_1 u_2 \right] \mbox{ dx }= (\mathbf{x}_1, \calG ^*\mathbf{x}_2)_{U \times L^2},
\end{align}
resulting in 
\begin{equation}
\label{gsgis}
\calG ^* = \left( \begin{array}{cc} \calA^{-1}  & \calA ^{-2} \\ -I & 0 \end{array} \right), \quad \calG^* \calG = \left( \begin{array}{cc} I + \calA ^{-2} & -\calA^{-1}  \\ -\calA  & I \end{array} \right).
\end{equation}

\subsection{Backwards heat problem.} \label{sec:bh}
The same analysis can be done for the backwards heat problem, although it is a bit trickier, due to time dependence.
The operator of interest in the all-at-once formulation of the backwards heat equation is given by
\begin{equation}
\calG  = \left( \begin{array}{cc} \partial_t +\calA & \underline{I}\calA \\ \delta_T & I \end{array} \right),
\end{equation}
which is an operator from the space $U_0\times X$ to $W'\times Y:=L^2(H^{-1}(\Omega) \times L^2(\Omega)$ where $X=H_0^1(\Omega)$,
\[
U_0 := \{ w \in L^2(H_0^1(\Omega)) \cap H^1(H^{-1}(\Omega)): w(0,x) = 0 \} \,,
\] 
which means, that the initial condition $(u+\ulI\theta)(0)=\theta$ is implicitly enforced through the function space we are using for $u$. The operator $\underline{I}\calA$ maps the static parameter into a time-dependent function space, i.e. $\theta \mapsto ( t \mapsto \calA \theta)$, and $I$ here denotes the embedding of $H_0^1(\Omega)$ into $L^2(\Omega)$. The scalar product on $U_0$ for all $v(t,x),u(t,x) \in U_0$ is given by
\begin{align}\label{innprodU0}
\langle v,u\rangle_{U_0} &=\int_0^T \int_{\Omega} \left[ \nabla v \nabla u + \nabla \calA ^{-1} v_t \nabla \calA ^{-1}  u_t \right] \mbox{dx dt} + \int_{\Omega} v(T) u(T) \mbox{ dx}  \\
 &=\int_0^T \int_{\Omega} \calA^{-1/2}(\partial_t+\calA)v \, \calA^{-1/2}(\partial_t+\calA) u  \mbox{ dx dt}  \\
&= \int_0^T \int_{\Omega} v \left[ \calA  u - \calA ^{-1} u_{tt}\right] \mbox{ dx dt} +\int_{\Omega} v(T) \left[u(T) + \calA ^{-1}u_t(T)\right] \mbox{ dx}.
\end{align}
The other scalar products of interest are those of $L^2(H^{-1})$ and of $H_0^1(\Omega)$:
\begin{equation}
\langle v, u\rangle_{L^2(H^{-1})} = \int_0^T \int_{\Omega} \nabla \calA ^{-1} v \cdot\nabla \calA ^{-1} u \mbox{ dx dt}\,, \quad 
\langle v, u\rangle_{H_0^1} = \int_{\Omega} \nabla v \cdot\nabla u \mbox{ dx}
\end{equation}
With the help of these scalar products the operator $\calG ^*$ will form as
\begin{equation*}
\calG ^* = \left( \begin{array}{cc} (\partial_t +\calA )^* & \delta_T^* \\ (\ulI\calA)^* & I^* \end{array} \right),
\end{equation*}
where for some components of the operator a bit more investigation is needed. We start with computing the first operator $(\partial_t +\calA )^*$ with the help of the scalar product equation
\begin{equation}
\label{adj1}
\langle(\partial_t +\calA )u, v\rangle_{L^2(H^{-1})} = \langle u, (\partial_t +\calA )^*v\rangle_{U_0},
\end{equation}
where the left hand side computes as
\begin{align*}
 \langle(\partial_t +\calA )u, v\rangle_{L^2(H^{-1})} &= \int_0^T \int_{\Omega} u(t,x) \left[ \calA ^{-1}v_t + v\right](t,x) \mbox{ dx dt} \\
&+ \int_{\Omega} u(T,x) \calA ^{-1}v(T,x) \mbox{ dx}
\end{align*}
and the right hand side, with $(\partial_t +\calA )^*v =: z$ according to the identity \eqref{innprodU0} works out as
\begin{align*}
\langle u, (\partial_t +\calA )^*v\rangle_{U_0} = \int_0^T \int_{\Omega} u(t,x) \left[ \calA  z - \calA ^{-1} z_{tt}\right](t,x) \mbox{ dx dt} \\
+\int_{\Omega} u(T,x) [z(T,x) + \calA ^{-1}z_t(T,x)] \mbox{ dx}. 
\end{align*}
Then, equation (\ref{adj1}) leads to the system
\begin{equation}\label{sys1}
\begin{array}{rcll}
\calA ^2 z -  z_{tt} &=& \calA v -v_t  \quad &\mbox{in } (0,T)\\
\calA  z(T) + z_t(T)  &=& v (T)  
\end{array}
\end{equation}
By factorizing $\calA ^2- \partial_{tt} = (\calA  -\partial_t)(\calA  + \partial_t)$	 and from the fact that $z \in U_0$ the system in (\ref{sys1}) results in
\begin{equation*}
\begin{array}{rcll}
(\partial_t+\calA)z &=& v \quad &\mbox{in } (0,T) \\
z(0) &=& 0 \,.
\end{array}
\end{equation*}
Therefore the adjoint operator $(\partial_t +\calA )^*$ applied to $v$ can be written as the solution of the heat equation using the variation of constants formula for the heat semigroup
\begin{equation*}
((\partial_t +\calA )^*v)(t,x) = \int_0^t e^{-\calA (t-s)}v(s,x) \mbox{ ds}.
\end{equation*}
The next operator computed is $\delta_T^*$. We proceed as before with the scalar product equation
\begin{equation*}
(\delta_T u, v)_{L^2(\Omega)} = (u, \delta_T^* v)_{U_0},
\end{equation*}
which, with rewriting $z = \delta_T^*v(x)$, due to \eqref{innprodU0} is equivalent to
\begin{align}
\label{deltats}
\int_{\Omega} u(T,x)v(x) \mbox{ dx} = \int_0^T \int_{\Omega} u(t,x)[\calA  z - \calA ^{-1}z_{tt}](t,x) \mbox{ dx dt} \\
+ \int_{\Omega} u(T,x)[z(T,x) + \calA^{-1}z_t(T,x)] \mbox{ dx}.
\end{align}
Then, (\ref{deltats}) leads to the system 
\begin{equation*}
\begin{array}{rcll}
(\partial_t-\calA  )(\partial_t+\calA)z &=& 0 \quad &\mbox{in } (0,T) \\
(\partial_t+\calA)z(T) &=& \calA  v \,.
\end{array}
\end{equation*}
Again using the heat semigroup, the solution can be given as
\begin{equation*}
(\delta_T^* v)(t,x) = \calA  \int_0^t e^{-\calA (T+t-2s)}v(x) \mbox{ ds} 
= \tfrac12[e^{-\calA (T-t)} - e^{-\calA (T+t)}]v(x).
\end{equation*}
The other two adjoint operators compute as
\begin{equation*}
((\calA\ulI)^*u)(x) = \calA ^{-1} \int_0^Tu(s,x) \mbox{ ds}, \quad (I^*v)(x) = (\calA ^{-1}v)(x).
\end{equation*}
Altogether, the adjoint operator $\calG ^*$ is given by
\begin{equation*}
\calG ^* = \left( \begin{array}{cc} \int_0^\cdot e^{-\calA  (\cdot-s)} .(s) \mbox{ ds} & \tfrac12[e^{-\calA (T-\cdot)} - e^{-\calA (T+\cdot)}]\\ \calA ^{-1}\ulI^* & \calA^{-1} \end{array} \right),
\end{equation*}
and $\calG ^* \calG $ is given by
\begin{equation*}
\calG ^* \calG  = \left( \begin{array}{cc} 
I +  \tfrac12[e^{-\calA (T-\cdot)} - e^{-\calA (T+\cdot)}] \delta_T \quad& 
\ulI-e^{-\calA \cdot} + \tfrac12[e^{-\calA (T-\cdot)} - e^{-\calA (T+\cdot)}] \\ 
\ulI^* + 2 \calA ^{-1} \delta_T & 
TI + \calA ^{-1} \end{array} \right).
\end{equation*}

\section{Analysis of the eigenvalues.} \label{sec:eig}
In this section the eigenvalues of the all-at-once operators are analyzed in terms of the two prototypical examples. For the analysis, the adjoint operators from Section \ref{sec:adj} will be used. Especially the operator $\calG ^*\calG $ is of interest. As one might have already noticed, these operators do not look symmetric. However, they are indeed symmetric, but with respect to the specific inner product on $U\times L^2(\Omega)$ and $U_0\times H_0^1(\Omega)$, respectively. Therefore, a transformation is applied first, to find a representation of the operator in an $L^2$-related inner product. It can be shown that this transformed operator leads to approximate eigenvalues of the true operator.\\

\begin{lemma}\label{lem:eigvaltransf}
Let $V$ and $H$ be Hilbert spaces with $C:V \to V$ self-adjoint and compact and $\mathcal{T}:V \to H$ boundedly invertible with $\mathcal{T}\in L(V,H)$ and $\mathcal{T}^{-1} \in L(H,V)$. Then the operator $\tilde{C}:= (\mathcal{T}^{-1})^*C \mathcal{T}^{-1} : H \to H$ is self-adjoint and compact and the eigenvalues $\lambda_k$ of $C$ and $\mu_k$ of $\tilde{C}$ decay at the same rate, more precisely it holds
\begin{equation*}
\frac{1}{\|\mathcal{T}^{-1}||^2}\mu_k \leq \lambda_k \leq \|\mathcal{T}\|^2 \mu_k,
\end{equation*}
with $\lambda_1 \geq \lambda_2 \geq \cdots \geq 0, \mu_1 \geq \mu_2 \geq \cdots \geq 0$.
\end{lemma}

    {\em Proof.} \
The proof is based on the Courant-Fischer Theorem, which we quote here for the convenience of the reader
\begin{quote}
\textsc{Theorem (Courant-Fischer).} Let $C:V\to V$ be a selfadjoint and compact operator. Then the eigenvalues of $C$ fulfill
\begin{equation*}
\lambda_k = \max \{ \min\{(Cx,x)_V: x\in S_k, \|x\|=1\}: \dim (S_k)=k, S_k \mbox{ subspace of V}\} 
\end{equation*}
\end{quote}
 Let
\begin{align*}
\lambda_k &= \max \{ \min \{ (Cx,x)_V : x \in S_k, \|x\|=1 \}: \dim (S_k)=k, S_k \mbox{ subspace of }V \} = \\
&= \max_{\dim(S_k)=k} \min_{x \in S_k, \|x\| = 1} (Cx,x)_V =  \max_{\dim(S_k)=k} \min_{x \in S_k, \|x\|=1} ((\mathcal{T}^{-1})^*C\mathcal{T}^{-1}\, \mathcal{T}x, \mathcal{T}x)_V = \\
&=\max_{\dim(S_k)=k} \min_{x \in S_k, \hat{x} = \mathcal{T}x / \|\mathcal{T}x\|, \|x\|=1} (\tilde{C}\hat{x}, \hat{x})_H \|\mathcal{T}x\|^2 \quad (\star).
\end{align*}
Due to the fact that
\begin{equation*}
\hat{x} = \frac{\mathcal{T}x}{\|\mathcal{T}x\|} \in \hat{S}_k = \mathcal{T}S_k,
\end{equation*}
and the dimension of $S_k$ being $k$, due to regularity of $\mathcal{T}$, $\hat{S}_k$ is of dimension $k$ as well. Therefore, taking the minimum over a superset and using $\|\mathcal{T}x\|\geq\frac{1}{\|\mathcal{T}^{-1}\|}\|x\|$ results in
\begin{equation*}
 (\star) \geq \frac{1}{\|\mathcal{T}^{-1}\|^2} \max_{\dim(\hat{S}_k)=k} \min_{ \hat{x} \in \hat{S}_k, \|x\| =1}(\tilde{C}x,x)_H = \frac{1}{\|\mathcal{T}^{-1}\|^2} \mu_k.
\end{equation*}
Analogously it holds
\begin{equation*}
\mu_k  \geq \frac{1}{\|(\mathcal{T}^{-1})^{-1}\|^2} \lambda_k = \frac{1}{\|\mathcal{T}\|^2} \lambda_k.
\end{equation*}
    $\diamondsuit$\\

With the help of Lemma~\ref{lem:eigvaltransf} we will transform the operator $\calG ^*\calG $ both for the inverse source and backwards heat problem and then investigate in the computation of the eigenvalues of the resulting operators.

\subsection{Inverse source problem.} The operator $\calG ^*\calG $ as stated in (\ref{gsgis}), will be transformed according to Lemma~\ref{lem:eigvaltransf} with the operator $\mathcal{T}:U\times L^2(\Omega) \to L^2(\Omega) \times L^2(\Omega)$ given by
\begin{equation*}
\mathcal{T} = \left( \begin{array}{cc} \calA  & 0 \\ 0 & I \end{array} \right),
\end{equation*}
which can be easily seen to be unitary in this setting of spaces, i.e., $(\mathcal{T}^{-1})^*=\mathcal{T}$,
therefore $\calG ^*\calG $ transforms to
\begin{eqnarray*}
\widetilde{\calG ^*\calG } = (\mathcal{T}^{-1})^*\calG ^*\calG  \mathcal{T}^{-1} = \left( \begin{array}{cc} \calA  & 0 \\ 0 & I \end{array} \right) \left( \begin{array}{cc} I + \calA ^{-2} & -\calA ^{-1}  \\ -\calA  & I \end{array} \right)  \left( \begin{array}{cc} \calA^{-1}  & 0 \\ 0 & I \end{array} \right) \\
= \left( \begin{array}{cc} I+ \calA ^{-2} & -I \\ -I & I \end{array} \right),
\end{eqnarray*}
which is obviously a compact perturbation of 
$\overline{\calG ^*\calG }=\left( \begin{array}{cc} I& -I \\ -I & I \end{array} \right)$ whose eigenvalues are $0$ and $2$ with eigenspaces 
$E_0=\{(g,g)^T\, : \, g\in L^2(\Omega)\}$,
$E_2=\{(-g,g)^T\, : \, g\in L^2(\Omega)\}$, 
that actually span all of $L^2(\Omega)$, $E_0\oplus E_2 = L^2(\Omega)$.

\subsubsection{Analytic computation of the eigenvalues.}
The eigenvalues $\lambda$ are computed both analytically and numerically. For the analytic computation the characteristic eigenvalue equation is considered
\begin{eqnarray*}
\left( \begin{array}{cc} I + \calA ^{-2} & -I \\ -I & I  \end{array} \right) \left( \begin{array}{c} f \\ g \end{array} \right) = \lambda \left( \begin{array}{c}  f \\  g\end{array} \right),
\end{eqnarray*}
which leads to the system
\begin{eqnarray*}
(1- \lambda) f + \calA ^{-2}f -g = 0, \quad (1-\lambda)g -f = 0,
\end{eqnarray*}
resulting in
\begin{eqnarray*}
\calA ^{-2}g = \frac{\lambda(2-\lambda)}{(1-\lambda)}g \\
f = (1-\lambda)g.
\end{eqnarray*}
Therefore, let $\mu = \lambda \frac{2-\lambda}{1 - \lambda}$ be an eigenvalue of $\calA ^{-2}$. Then $\lambda$ solves the quadratic equation
\begin{equation*}
\lambda^2 -\lambda(2+\mu) +\mu = 0,
\end{equation*}
resulting in the two solutions
\begin{eqnarray*}
\lambda_{1,2} = \frac{2+\mu}{2} \pm \sqrt{\frac{\mu^2}{4} +1} 
= 1+\frac{\mu}{2}+\sqrt{1+\frac{\mu^2}{4}} 
\,, \ 
\frac{\mu}{ 1+\frac{\mu}{2}+\sqrt{1+\frac{\mu^2}{4}} }.
\end{eqnarray*}
As $\mu \to 0$, the solutions for $\lambda_{1,2}$ tend to $2$ and $0$, respectively. Therefore, the eigenvalues of $\widetilde{\calG ^* \calG }$ accumulate at 0 and 2, at the same (basically linear) speed as $\mu \to 0$.

Thus, besides the singular values tending to zero at a linear rate, known from the reduced setting and reflecting the mild ill-posedness of the inverse source problem, we have another sequence tending to a positive value.

\subsubsection{Numerical computation of the eigenvalues.}
\begin{figure}
\centering
\begin{subfigure}[c]{0.45\textwidth}
\includegraphics[width=0.9\textwidth]{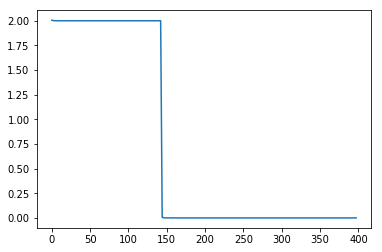}
\subcaption{Eigenvalues}
\label{fig:eig_is}
\end{subfigure}
\begin{subfigure}[c]{0.45\textwidth}
\includegraphics[width=0.9\textwidth]{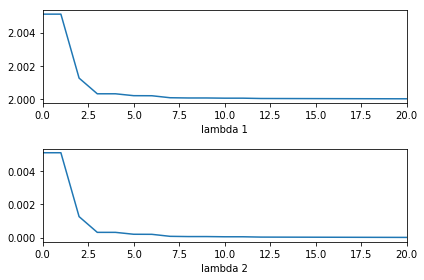}
\subcaption{Detailed plot}
\label{fig:rooteig_is}
\end{subfigure}
\caption{Numerical results for the eigenvalues of $\widetilde{\mathcal{G}^*\mathcal{G}}$ for the inverse source problem.}
\end{figure}

The computation is done in python with the help of the finite element discretization in FEniCS and the eigensolver from SLEPc. The discretzation is done on a unit square mesh with degree 1 Lagrange elements. To overcome the computation of $\calA ^{-2}$ the inverse operator 
$\widetilde{\calG ^* \calG }^{-1} = \left( \begin{array}{cc} \calA^2 & \calA^2 \\ \calA^2 & \calA^2 + I \end{array} \right)$ is used. 
The matrices resulting from the finite element discretization 
$u(x)\approx \sum_{i=1}^{n_{el}} u_i \phi_i(x)$ 
with the FE basis functions $\phi_i$, and the coefficient vector $\underline{u}=(u_1,\ldots,u_{n_{el}})$, according to the identities
\[
\calA u = \lambda u\ \Leftrightarrow\  K\underline{u}=\lambda M\underline{u}\,, \quad 
\calA^2 u = \lambda u\ \Leftrightarrow\  KM^{-1}K\underline{u}=\lambda M\underline{u}\,, 
\]
are
\begin{equation*}
\widetilde{\calG ^* \calG }^{-1}_h = \left( \begin{array}{cc} KM^{-1}K & KM^{-1}K \\ KM^{-1}K & KM^{-1}K + M \end{array} \right), \quad 
\mathcal{M} = \left( \begin{array}{cc} M & 0 \\ 0 & M \end{array} \right),
\end{equation*}
with $M$ and $K$ being the assembled mass and stiffness matrices according to 
\begin{eqnarray}
M_{i,j}&=&\int_\Omega \phi_i\phi_j \label{M},\\ 
K_{i,j}&=&\int_\Omega \nabla\phi_i\cdot\nabla\phi_j \label{K}, 
\end{eqnarray}
to state the eigenvalue equation as
\begin{equation*}
\mathcal{M} \mathbf{x} = \lambda \widetilde{\calG ^* \calG }^{-1}_h \mathbf{x},
\end{equation*}
with $\lambda$ and $\mathbf{x}$ denoting the eigenvalue and eigenvector, respectively.
Sampling of 500 eigenvalues leads then to the following visualized output as seen in Figure \ref{fig:eig_is}, the square root scaled plot in Figure \ref{fig:rooteig_is} suggests quadratic decay, corresponding to the fact that $\lambda_{2,n}\sim\frac{\mu_n}{2}\sim\frac{C}{n^2}$ as the eigenvalues of the Laplacian on the 2-d unit square  are given by $\mu_{j,k}=\frac{j^2+k^2}{4}$, $j,k\in\mathbb{N}$ and therefore, upon proper renumbering, decay linearly in $n$.

\subsection{Backwards heat equation.} The same analysis is done with the backwards heat equation. Here, the operator of interest is
\begin{equation}
 \calG ^* \calG  = \left( \begin{array}{cc} I +  \calA \int_0^t e^{-\calA (T+t-2s)} \mbox{ .(s)} \delta_T & I-e^{-\calA t} + \calA \int_0^t e^{-\calA (T+t-2s)} \mbox{ .(s)} \\ \ulI^* + 2 \calA ^{-1} \delta_T & TI + \calA ^{-1} \end{array} \right).
\end{equation}
The transformation from the space $U_0 \times H_0^1(\Omega)$ to the space $L^2(L^2(\Omega)) \times L^2(\Omega)$ is computed with Lemma~\ref{lem:eigvaltransf} where the operator $\mathcal{T}$ for the transformation is given by
\begin{equation*}
\mathcal{T} = \left( \begin{array}{cc} \calA ^{-1/2}(\partial_t + \calA ) & 0 \\ 0 & \calA ^{1/2} \end{array} \right)
\end{equation*}
and therefore,
\[
\widetilde{\calG ^* \calG } = (\calG  \mathcal{T}^{-1})^*\calG  \mathcal{T}^{-1}.
\]
Here  
\[
\mathcal{T}^{-1} \left( \begin{array}{c} f \\ g \end{array} \right) = \left( \begin{array}{c} u \\ \theta \end{array} \right) \ \Leftrightarrow \
\left\{ \begin{array}{c} (\partial_t + \calA )u = \calA ^{1/2} f \\ \theta = \calA ^{-1/2} g \end{array} \right. 
\]
and
\[
\calG  \left( \begin{array}{c} u \\ \theta \end{array} \right) = \left( \begin{array}{c} a \\ b \end{array} \right) \ \Leftrightarrow \
\left\{ \begin{array}{c} (\partial_t + \calA )u + \calA  \theta = a \\ u(T)+\theta = b \end{array} \right., 
\]
thus
\[
\left\{ \begin{array}{c} \calA ^{1/2} f(t) + \calA ^{1/2} g = a(t), \ t\in(0,T) \\ 
\int_0^T e^{-\calA (T-t)}\calA ^{1/2} f(t)\, dt + \calA ^{-1/2} g = b
\end{array} \right.\,,
\]
i.e., 
\[
\calG \mathcal{T}^{-1} = \left( \begin{array}{cc} \calA ^{1/2} & \calA^{1/2}\ulI\\
\int_0^T e^{-\calA (T-s)} \calA ^{1/2} .(s) \mbox{ ds} & \calA ^{-1/2}
\end{array} \right).
\]
For computing the adjoint $(\calG \mathcal{T}^{-1})^*$ we consider the identity
\begin{align*}
& \left\langle
\calG \mathcal{T}^{-1} \left( \begin{array}{c} f\\g \end{array} \right),
\left( \begin{array}{c} a\\b \end{array} \right)
\right\rangle_{L^2(H^{-1})\times L^2} &\\
=&\int_0^T\int_\Omega (f(t,x)+g(x))\calA^{-1/2}a(t,x)\mbox{ dx}\mbox{ dt} \\
&+\int_\Omega \left(\int_0^T e^{-\calA(T-s)}\calA^{1/2}f(s,x)\mbox{ ds} +\calA^{-1/2}g(x)\right)b(x)\mbox{ dx}\\
=&\int_0^T\int_\Omega f(t,x)\left(\calA^{-1/2}a(t,x)+e^{-\calA(T-t)}\calA^{1/2}b\right)\mbox{ dx}\mbox{ dt} \\
&+\int_\Omega g(x)\calA^{-1/2}\left(\int_0^T a(t,x)\mbox{ dt} + b\right)\mbox{ dx}, 
\end{align*}
which yields
\[
(\calG \mathcal{T}^{-1})^*=\left( \begin{array}{cc} \calA ^{-1/2} & e^{-\calA (T-\cdot)} \calA^{1/2}\\ \calA ^{-1/2}\ulI^* & \calA ^{-1/2}
\end{array} \right).
\]
Therefore,
\begin{align}
& \widetilde{\calG ^* \calG } = (\calG \mathcal{T}^{-1})^*  \calG \mathcal{T}^{-1}  = \nonumber\\
 &= \left( \begin{array}{cc} I + e^{-\calA (T-\cdot)} \calA ^{1/2} \int_0^T e^{-\calA (T-s)} \calA ^{1/2} .(s) \mbox{ ds} & \ulI + e^{-\calA (T-\cdot)} \\ \ulI^* + \int_0^T e^{-\calA (T-s)} .(s) \mbox{ ds} & TI + \calA ^{-1} \end{array} \right).
\label{bhtilGG}
\end{align}
This shows that $\widetilde{\calG ^* \calG }$ is a compact perturbation of the operator 
$\overline{\calG ^* \calG }=\left( \begin{array}{cc} I+D&\ulI\\\ulI^*& TI \end{array} \right)$,
where $D= e^{-\calA (T-\cdot)} \calA ^{1/2} \int_0^T e^{-\calA (T-s)} \calA ^{1/2} .(s) \mbox{ ds}$.
To see this, consider a system of eigenvalues $\gamma_n$ and eigenfunctions $\phi_n$ of $\calA$ and the estimate
\begin{eqnarray*}
&&\|\int_0^T e^{-\calA (T-s)} u(s) \mbox{ ds}\|_{\dot{H}^1(\Omega)} =
\|\calA^{1/2}\int_0^T e^{-\calA (T-s)} u(s) \mbox{ ds}\|_{\dot{H}^1(\Omega)} \\
&&= \left(\sum_{n=1}^\infty \gamma_n \Bigl(\int_0^T e^{-\gamma_n (T-s)} \langle u(s),\phi_n\rangle_{L^2}\mbox{ ds} \Bigr)^2\right)^{1/2}\\
&&\leq \Bigl(\sum_{n=1}^\infty \int_0^T \langle u(s),\phi_n\rangle_{L^2}^2 \mbox{ ds}  
\underbrace{\gamma_n \int_0^T e^{-2\gamma_n (T-s)} \mbox{ ds}}_{=\frac12(1-e^{-2\gamma_nT})} 
\Bigr)^{1/2}
\ \leq \frac{1}{\sqrt{2}} \|u\|_{L^2(L^2)}\,,
\end{eqnarray*}
which together with compactness of the embedding $\dot{H}^1(\Omega):=\mathcal{D}(\calA^{1/2})\to L^2(\Omega)$, shows that the operator  $\int_0^T e^{-\calA (T-s)} .(s):L^2(0,T;L^2(\Omega))\to L^2(\Omega)$ is compact and so is its adjoint $e^{-\calA (T-\cdot)}:L^2(\Omega)\to L^2(0,T;L^2(\Omega))$. However $D:L^2(0,T;L^2(\Omega))\to L^2(0,T;L^2(\Omega))$ is only bounded, as the following computation shows
\begin{eqnarray*}
&&\|e^{-\calA (T-\cdot)} \calA ^{1/2} \int_0^T e^{-\calA (T-s)} \calA ^{1/2} u(s) \mbox{ ds}\|_{L^2} \\
&&= \left(\int_0^T\sum_{n=1}^\infty \Bigl(\gamma_n e^{-\gamma_n (T-t)} \int_0^T e^{-\gamma_n (T-s)} \langle u(s),\phi_n\rangle_{L^2} \mbox{ ds}\Bigr)^2 \mbox{ dt}\right)^{1/2}\\
&&\leq \left(\int_0^T\sum_{n=1}^\infty \gamma_n^2 e^{-2\gamma_n (T-t)} \int_0^T e^{-2\gamma_n (T-s)}\, ds \ \int_0^T \langle u(s),\phi_n\rangle_{L^2}^2 \mbox{ ds}  \mbox{ dt}\right)^{1/2}\\
&& = \left(\sum_{n=1}^\infty \int_0^T\langle u(s),\phi_n\rangle_{L^2}^2 \mbox{ ds} \Bigl(\gamma_n e^{-2\gamma_n (T-t)} \mbox{ dt}\Bigr)^2 \right)^{1/2}
\ \leq \frac12 \|u\|_{L^2(L^2)}\,.
\end{eqnarray*}

\subsubsection{Analytic computation of the eigenvalues.}
With that operator we now again state the eigenvalue equation
\begin{eqnarray}
\label{bhev}
\widetilde{\calG ^*\calG } \left( \begin{array}{c} f \\ g \end{array} \right) = \lambda \left( \begin{array}{c} f \\ g \end{array} \right),
\end{eqnarray}
where we write $f(t)$ and $g$ in terms of their generalized Fourier series with respect to the eigensystem $(\mu_n,\phi_n)_{n\in\mathbb{N}}\subseteq\mathbb{R}\times L^2(\Omega)$ of $\calA^{-1}$
\[
f(t,x)=\sum_{n=1}^\infty f_n(t)\phi_n(x)\,,\quad 
g(x)=\sum_{n=1}^\infty g_n\phi_n(x),
\]
with $f_n(t)=\langle f(t),\phi_n\rangle_{L^2(\Omega)}$, $g_n=\langle g,\phi_n\rangle_{L^2(\Omega)}$. 
Then, with \eqref{bhtilGG}, upon taking inner products with $\phi_n$, equation (\ref{bhev}) leads to the system
\begin{equation}
\label{eig_fngn}
\begin{array}{l}
(1-\lambda) f_n(t) + e^{-\frac{1}{\mu_n}(T-t)} \frac{1}{\mu_n} \int_0^T e^{-\frac{1}{\mu_n}(T-s)} f_n(s) \mbox{ ds} + (1+e^{-\frac{1}{\mu_n}(T-t)})g_n  = 0, \\
\int_0^T f_n(s) \mbox{ ds} + \int_0^Te^{-\frac{1}{\mu_n}(T-s)} f_n(s) \mbox{ ds} + (T+\mu_n-\lambda)g_n = 0,
\end{array}
\end{equation}
for all $n\in\mathbb{N}$ and $t \in (0,T)$.
This has a nontrivial solution $(f_n,g_n)_{n\in\mathbb{N}}\in L^2(0,T,\ell^2)\times\ell^2$, if and only if there exists $m\in\mathbb{N}$ such that \eqref{eig_fngn} holds for $n=m$ and $f_m\not\equiv0$ or $g_m\not=0$ (to see sufficiency, set all other components of $(f_n,g_n)_{n\in\mathbb{N}}$ to zero and they will trivially satisfy the linear system \eqref{eig_fngn}).
With the particular time dependent functions $b^m_0:t\mapsto 1$ and $b^m_1:t\mapsto e^{-\frac{1}{\mu_n}(T-t)}$ playing a role here, this is the case iff either 
\[
\mbox{(a) } \lambda=1 \mbox{ and }
\left\{
\begin{array}{l}
g_m b^m_0(t)+\Bigl(\tfrac{1}{\mu_m}\int_0^T b^m_1(s) f_m(s)\mbox{ ds}+g_m\Bigr) b^m_1(t) = 0 \,, \ t\in(0,T)\\
\int_0^T b^m_0(s)f_m(s)\mbox{ ds} + \int_0^T b^m_1(s) f_m(s)\mbox{ ds}+(T+\mu_m-\lambda)g_m = 0
\end{array}\right.
\]
or 
\[
\mbox{(b) } \lambda\not=1 \mbox{ and }
\left\{
\begin{array}{l}
f_m = a_m b^m_0 + c_m b^m_1 \mbox{ where}\\
a_m=\tfrac{1}{\lambda-1}g_m \,, \quad c_m =\tfrac{1}{\lambda-1}\Bigl(\tfrac{1}{\mu_m}\int_0^T b^m_1(s) f_m(s)\mbox{ ds}+g_m\Bigr)\\
\int_0^T b^m_0(s)f_m(s)\mbox{ ds} + \int_0^T b^m_1(s) f_m(s)\mbox{ ds}+(T+\mu_m-\lambda)g_m = 0
\end{array}\right..
\]
The first equation of case (a) $\lambda=1$ due to linear independence of the functions $b^m_0$, $b^m_1$ leads to the conditions $g_m=0$ and $\tfrac{1}{\mu_m}\int_0^T b^m_1(s) f_m(s)\mbox{ ds}+g_m=0$,  that, combined with the second equation, yield the necessary and sufficient conditions
\[
g_m=0, \int_0^T b^m_1(s) f_m(s)\mbox{ ds}=0 \mbox{ and } \int_0^T b^m_0(s) f_m(s)\mbox{ ds}=0\,.
\]
This corresponds to a unit eigenvalue with the infinite dimensional eigenspace $E_1$ 
\begin{equation}\label{lambda1}
\lambda_0=1\,, \quad
E_1 = \mbox{span}\{(x,t)\mapsto \psi(t)\phi_m(x)\ : \ \psi\in\{b^m_0,\, b^m_1\}^\bot\,, \ m\in\mathbb{N}\}
\end{equation}
where ${.}^\bot$ denotes the $L^2(0,T)$ orthogonal complement.
\\
In case (b) $\lambda_m\not=1$ inserting the representation $f_m = a_m b^m_0 + c_m b^m_1$ into the two equations involving $f_m$ and $g_m=(\lambda-1)a_m$, we get the following system for $a_m$ and $b_m$
\begin{eqnarray*}
&& c_m =\tfrac{1}{\lambda-1}\tfrac{1}{\mu_m}\int_0^T b^m_1(s) (a_mb^m_0(s)+c_mb^m_1(s))\mbox{ ds}+a_m,\\
&& \int_0^T b^m_0(s)(a_mb^m_0(s)+c_mb^m_1(s))\mbox{ ds} + \int_0^T b^m_1(s) (a_mb^m_0(s)+c_mb^m_1(s))\mbox{ ds} \\
&& \quad +(T+\mu_m-\lambda)(\lambda-1)a_m = 0\,,
\end{eqnarray*}
which with the integrals $\int_0^T b^m_0(s)^2\mbox{ ds}=T$, $\int_0^T b^m_1(s)^2\mbox{ ds}=\frac{\mu_m}{2}(1-e^{-\frac{1}{\mu_m}2T})$, $\int_0^T b^m_0(s) b^m_1(s)\mbox{ ds}=\mu_m(1-e^{-\frac{1}{\mu_m}T})$ reads as
\begin{eqnarray*}
&&\Bigl((\lambda-1)+1-e^{-\frac{1}{\mu_m}T}\Bigr) a_m 
+ \Bigl(-(\lambda-1)+\tfrac12(1-e^{-\frac{1}{\mu_m}2T})\Bigr)c_m=0,\\
&& \Bigl(T+(T+\mu_m-\lambda)(\lambda-1)+\mu_m(1-e^{-\frac{1}{\mu_m}T})\Bigr) a_m \\
&& \quad + \mu_m(\tfrac32-e^{-\frac{1}{\mu_m}T}-e^{-\frac{1}{\mu_m}2T}) c_m = 0\,,
\end{eqnarray*}
i.e., 
\begin{eqnarray*}
&&(\lambda-e^{-\frac{1}{\mu_m}T}) a_m 
+(\tfrac32-\lambda-\tfrac12 e^{-\frac{1}{\mu_m}2T}))c_m=0,\\
&& \Bigl(-\lambda^2 +\lambda (T+1+\mu_m) -\mu_m e^{-\frac{1}{\mu_m}T}\Bigr) a_m
+ \mu_m(\tfrac32-e^{-\frac{1}{\mu_m}T}-e^{-\frac{1}{\mu_m}2T}) c_m = 0\,.
\end{eqnarray*}
Existence of a nontrivial solution $a_m$, $c_m$ by setting the determinant of this system to zero is equivalent to the following cubic equation for $\lambda$
\[
\lambda^3-\lambda^2(T+\tfrac52+\alpha_m)
+\lambda\left(\tfrac32(T+1)+\beta_m\right)
=\mu_me^{-\frac{1}{\mu_m}2T} 
\,,
\]
with $\alpha_m = \mu_m-\tfrac12 e^{-\frac{1}{\mu_m}T}$, $\beta_m = 2\mu_me^{-\frac{1}{\mu_m}T}-\tfrac12 e^{-\frac{1}{\mu_m}2T}(T+1)$
whose solutions $\lambda_{\mu_m}^1$, $\lambda_{\mu_m}^2$, $\lambda_{\mu_m}^3$ will be the remaining (besides $\lambda_0=1$) eigenvalues of $\widetilde{\mathcal{G}^*\mathcal{G}}$.
To investigate the asymptotics of the eigenvalues, in particular of $\lambda^1_{\mu_n}$, we consider another sequence of values, namely 
\begin{eqnarray*}
\bar{\lambda}_{1,m} &=&0=:\bar{\lambda}_1,\\
\bar{\lambda}_{2,m} &=& \tfrac{T+\tfrac52+\alpha_m}{2}+\tfrac{\sqrt{(T-\tfrac12)^2+
\alpha_m^2+(2T+5)\alpha_m-4\beta_m}}{2}
\to 
\left\{\begin{array}{l}T+1\mbox{ if }T\geq\tfrac12\\ \tfrac32\mbox{ if }T<\tfrac12\end{array}\right.
=:\bar{\lambda}_2\,, \\
\bar{\lambda}_{3,m} &=& \tfrac{\tfrac32(T+1)+\beta_m}{\bar{\lambda}_{2,m}}
\to 
\left\{\begin{array}{l}\tfrac32\mbox{ if }T\geq\tfrac12\\ T+1\mbox{ if }T<\tfrac12\end{array}\right.
=:\bar{\lambda}_3\,.
\end{eqnarray*}
(where the limits are to be understood as $m\to\infty$ and therewith $\mu_m\to0$),
so that the cubic equation for $\lambda$ above can be written as
\begin{equation}\label{lambda_asymp}
\lambda (\lambda - \bar{\lambda}_{2,m})(\lambda - \bar{\lambda}_{3,m})= \mu_m e^{-\frac{1}{\mu_m}2T}\,.
\end{equation}
Since $\lambda^i_{\mu_m}$ are roots of cubic polynomial whose coefficients converge to those of a cubic polynomial with the three single roots $\bar{\lambda}_1$, $\bar{\lambda}_2$, $\bar{\lambda}_3$, we also have the convergence 
\[
\lambda^i_{\mu_m}\to\bar{\lambda}_i\mbox{ as }m\to\infty\,, \ i\in\{1,2,3\}\,.
\] 
In particular this means that since $\bar{\lambda}_1=0<\bar{\lambda}_2,\bar{\lambda}_3$ and from \eqref{lambda_asymp} we have
\[
\lambda^1_{\mu_m} = \frac{\mu_m e^{-\frac{1}{\mu_m}2T}}{(\lambda^1_{\mu_m} - \bar{\lambda}_{2,m})(\lambda^1_{\mu_m} - \bar{\lambda}_{3,m})}\,,
\]
where the denominator on the right hand side is positive and bounded away from zero for $m$ sufficiently large. Hence, there exists $m_0\in\mathbb{N}$ such that for all $m\geq n_0$ we have $\lambda^1_{\mu_m}>0$ (compatibly with the fact that $\lambda^1_{\mu_m}$ is an eigenvalue of the positive semidefinite operator $\widetilde{\calG^*\calG}$) and $\lambda^1_{\mu_m} = O(\mu_m e^{-\frac{1}{\mu_m}2T})$.
\\
Thus, we here have, in addition to $\lambda_0=1$, three sequences of eigenvalues: One tending to zero at an exponential rate (like in the reduced setting) in according to the severe ill-posedness of the backwards diffusion problem, and two further sequences accumulating at the positive values $\tfrac32$ and $T+1$.

\subsubsection{Numerical computation of the eigenvalues.}
\begin{figure}[!tbp]
\centering
\includegraphics[width=.5\textwidth]{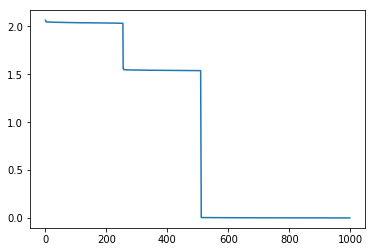}
\caption{Eigenvalues of $\widetilde{\calG^*\calG}$ for backwards heat problem.}
\label{figbh}
\end{figure}

The computation is done in python similarly to the computation of the eigenvalues of the inverse source problem, relying on a finite element discretization of the Laplacian. To obtain the right matrices for representing the semigroup expressions appearing in the definition of $\widetilde{\calG ^*\calG }$, we exemplarily consider the term $w(t)=\int_0^t e^{-(t-s)\calA}f(s)\, ds$ for some $f\in L^2(L^2(\Omega))$, and make a semidiscretization in space with a finite element ansatz $f(t,x)\approx\sum_{i=1}^{n_{el}} f_i(t)\phi_i(x)$, $w(t,x)\approx\sum_{i=1}^{n_{el}} w_i(t)\phi_i(x)$ with basis functions $\phi_1,\ldots,\phi_n$. Taking into account the fact that $w$ solves
\begin{equation*}
\begin{array}{rcll}
\partial_t w +\calA w &=& f \quad &\mbox{ in } (0,T)\times\Omega \\
w &=& 0 \quad &\mbox{ on } (0,T)\times\partial \Omega \\ 
w(0,x) &=& 0 \quad &x\in \Omega,
\end{array}
\end{equation*}
i.e., inserting the above ansatz, and testing with FE shape functions $\phi_j$, $j\in\{1,\ldots,n_{el}\}$,
(a procedure known as Faedo-Galerkin approximation,) 
we end up with the system of ODEs
\[ 
M\dot{\underline{w}}(t)+ K\underline{w}(t) = M \underline{f}(t)\,, t\in(0,T)\,, \quad 
\underline{w}(0)=0\,,
\]
where $M$ and $K$ are FE mass and stiffness matrices according to \eqref{M}, \eqref{K},
which, in order to obtain symmetry with respect to the Euclidean inner product, we write as
$M^{1/2}\dot{\underline{w}}(t)+ M^{-1/2}KM^{-1/2} \, M^{1/2}\underline{w}(t) = M^{1/2} \underline{f}(t)$.
Thus the coefficient vector function $\underline{w}(t)=(w_1(t),\ldots,w_{n_{el}}(t))^T$ is determined by the identity 
\[ 
M^{1/2} \underline{w}(t) = \int_0^t e^{-(t-s)A_h}M^{1/2}\underline{f}(s)\,, ds
\] 
with $A_h= M^{-1/2}KM^{-1/2}$.
\\
Additional discretization in time is done by piecewise constant basis functions on a uniform partition $0=t_0<t_1<\cdots<t_N=T$, $t_k=\frac{T}{N}j=\tau k$, of the time interval, i.e, $f_i(t)=\sum_{k=0}^N f_i^k\psi^k(t)$ with $\psi^0(t)=\chi_{[0,\tau/2)}$, $\psi^N(t)=\chi_{[T-\tau/2,T)}$, $\psi^k(t)=\chi_{[t_k-\tau/2,t_k+\tau/2)}$, $k=1,\ldots,N-1$, so that with
\[
f_h(t,x)=\sum_{k=0}^N \sum_{i=1}^{n_{el}} f_i^k \psi^k(t)\phi_i(x)\,, \quad
g_h(x)=\sum_{i=1}^{n_{el}} g_i \phi_i(x)\,,
\]
the discretized eigenvalue equation in variational form reads as 
\begin{align}\label{vareig}
\left\langle \widetilde{\calG^*\calG}\left(\begin{array}{c}f_h\\g_h\end{array}\right),
\left(\begin{array}{c}\psi^\ell\phi_j\\ \phi_j\end{array}\right)\right\rangle_{L^2(L^2)\times L^2}
= \left\langle \left(\begin{array}{c}f_h\\g_h\end{array}\right),\left(\begin{array}{c}\psi^\ell\phi_j\\ \phi_j\end{array}\right)\right\rangle_{L^2(L^2)\times L^2},\\
\qquad \forall \ell\in\{0,\ldots,N\}, j\in\{1,\ldots,n_{el}\},
\end{align}
with $\widetilde{\calG ^* \calG }$ as in \eqref{bhtilGG}, i.e., in matrix-vector form
\[
\left(\begin{array}{cc}\mathcal{M}+\mathcal{E}& (\mathcal{M}^1+\mathcal{E}^1)^T\\
\mathcal{M}^1+\mathcal{E}^1& T\,M+M^{1/2}A_h^{-1}M^{1/2}\end{array}\right)
\left(\begin{array}{c}\underline{\underline{f}} \\ \underline{g}\end{array}\right)
=\lambda\left(\begin{array}{cc}\mathcal{M}&0\\0&M\end{array}\right)
\left(\begin{array}{c}\underline{\underline{f}} \\\underline{g}\end{array}\right)\,,
\]
with the matrices
\begin{eqnarray*}
\mathcal{M}&=&\tau\mbox{diag}(\tfrac12 M,M,\cdots,M,\tfrac12 M)\,,\\
\mathcal{M}^1&=&\tau(\tfrac12 M,M,\cdots,M,\tfrac12 M)\,,\\
\mathcal{E}_{\ell,k}&=&\tau^2 M^{1/2} A_h e^{-A_h(2T-t_\ell-t_k)}M^{1/2}, \\
\mathcal{E}_{0,k} &=& \mathcal{E}_{k,0}=\tfrac{\tau^2}{2} M^{1/2}A_h e^{-A_h(2T-t_k)}M^{1/2}, \quad
\mathcal{E}_{0,0} =\tfrac{\tau^2}{4} M^{1/2}A_h e^{-A_h 2T}M^{1/2}, \\
\mathcal{E}_{N,k} &=& \mathcal{E}_{k,N}=\tfrac{\tau^2}{2} M^{1/2}A_h e^{-A_h(T-t_k)}M^{1/2}, \quad
\mathcal{E}_{N,N} =\tfrac{\tau^2}{4} M^{1/2}A_h M^{1/2}, \\
\mathcal{E}^1 &=& \tau (\tfrac12 M^{1/2}e^{-A_hT}M^{1/2},\ldots,M^{1/2}e^{-A_h(T-t_k)}M^{1/2},\ldots,\tfrac12 M)
\end{eqnarray*}
for $\ell,k\in\{1,\ldots,N-1\}$.
Here we have computed and approximated the time integrals in \eqref{vareig} as follows
\begin{align*}
&\int_{t_k-\tau/2}^{t_k+\tau/2} \mbox{ ds} =\tau, \\
&\int_{t_k-\tau/2}^{t_k+\tau/2} e^{-A_h(T-s)} \mbox{ ds} = A_h^{-1}\left(e^{A_h\tau/2}-e^{-A_h\tau/2}\right) e^{-A_h(T-t_k)}\approx \tau e^{-A_h(T-t_k)}, \\
&\int_{t_\ell-\tau/2}^{t_\ell+\tau/2} e^{-A_h(T-t)} A_h \int_{t_k-\tau/2}^{t_k+\tau/2} e^{-A_h(T-s)} \mbox{ ds}\mbox{ dt} =
A_h^{-1}\left(e^{A_h\tau}-2+e^{-A_h\tau}\right) e^{-A_h(2T-t_\ell-t_k)}\\
&\hspace*{6cm} \\
&\quad \approx \tau^2 A_h e^{-A_h(2T-t_\ell-t_k)}\,.
\end{align*}
To simplify the implementation we transform the eigenvalue equation by pre-multiplication with 
$\left(\begin{array}{cc} \mathcal{M} & 0 \\ 0 & M \end{array} \right)^{-\frac{1}{2}}$ 
and setting $\left( \begin{array}{c} \underline{\underline{\tilde{f}}} \\ \underline{\tilde{g}} \end{array} \right)
=\left(\begin{array}{cc} \mathcal{M} & 0 \\ 0 & M \end{array} \right)^{\frac{1}{2}}
\left( \begin{array}{c} \underline{\underline{f}} \\ \underline{g} \end{array} \right) $
which results in the following equation
\begin{eqnarray}
\left( \begin{array}{cc} \mathcal{I} + \tilde{\mathcal{E}} & (\mathcal{I}^1 + \tilde{\mathcal{E}}^1)^T \\ \mathcal{I}^1 + \tilde{\mathcal{E}}^1 & TI + A_h^{-1} \end{array} \right) \left( \begin{array}{c} \underline{\underline{\tilde{f}}} \\ \underline{\tilde{g}} \end{array} \right) = \lambda \left( \begin{array}{c} \underline{\underline{\tilde{f}}} \\ \underline{\tilde{g}} \end{array} \right) \,,
\end{eqnarray}
with the matrices
\begin{eqnarray*}
\mathcal{I}&=&\mbox{diag}(I,\cdots,I)\,,\\
\mathcal{I}^1&=&\sqrt{\tau}(\tfrac{1}{\sqrt{2}} I,I,\cdots,I,\tfrac{1}{\sqrt{2}} I)\,,\\
 \tilde{\mathcal{E}}_{\ell,k}&=& \tau A_h e^{-A_h(2T-t_\ell-t_k)}, \\
 \tilde{\mathcal{E}}_{0,k} &=& \tilde{\mathcal{E}}_{k,0}=\tfrac{\tau}{\sqrt{2}} A_h e^{-A_h(2T-t_k)}, \quad
\tilde{\mathcal{E}}_{0,0} =\tfrac{\tau}{2} A_h e^{-A_h 2T}, \\
\tilde{\mathcal{E}}_{N,k} &=& \tilde{\mathcal{E}}_{k,N}= \tfrac{\tau}{\sqrt{2}} A_h e^{-A_h(T-t_k)}, \quad
\tilde{\mathcal{E}}_{N,N} =\tfrac{\tau}{2} A_h, \\
 \tilde{\mathcal{E}}_{0,N}&=& \tilde{\mathcal{E}}_{N,0} = \tau A_h e^{-A_h}, \\
\tilde{\mathcal{E}}^1 &=& \sqrt{\tau} (\tfrac{1}{\sqrt{2}} e^{-A_hT},\ldots,e^{-A_h(T-t_k)},\ldots,\tfrac{1}{\sqrt{2}} I).
\end{eqnarray*}
With that setting the eigenvalues are computed with the time interval $[0,T]$ with $T=1$ discretized in the points $t_0=0, t_1=0.25, t_2=0.5, t_3=0.75, t_4=T=1$. 
Sampling of 700 eigenvalues leads to the following visualized output, which can be seen in Figure \ref{figbh}. The eigenvalues tend to zero, in three steps, which was also found in the analytic considerations.

\section{Convergence analysis.} \label{sec:conv}
For the convergence analysis, the results in \cite{agapioumathe2018, mathe2018} are extended to problems formulated in an all-at-once fashion, which can be done in a very straightforward fashion, so we keep this section short by more or less recalling the essential results from \cite{agapioumathe2018, mathe2018}. 

The aim is to find out about the convergence of the posterior towards the true element $\mathbf{x}^*$ which generates the data $\mathbf{y}^{\delta}$. This is typically done by analyzing the squared posterior contraction, given by
\begin{equation}
\mbox{SPC } := \mathbb{E}^{\mathbf{x}^*} \mathbb{E}_{\alpha}^{\delta} \|\mathbf{x}^* -\mathbf{x}\|, \quad \alpha, \delta >0,
\end{equation}
where the outward expectation is taken with respect to the data generating function, which means with given $\mathbf{x}^*$, the distribution which generates the data $\mathbf{y}^{\delta}$. The inward expectation is taken with respect to the posterior distribution with given data $\mathbf{y}^{\delta}$ and chosen scaling parameter $\alpha$. The squared posterior contraction can be decomposed into the squared bias, the estimation variance and the posterior spread, i.e.
\begin{equation}
\mbox{SPC } = 
\|\mathbf{x}^* - \mathbb{E}^{\mathbf{x}^*}\mathbf{x}_{\alpha}^{\delta} \|^2 
+ \mathbb{E}^{\mathbf{x}^*}\|\mathbf{x}_{\alpha}^{\delta} - \mathbb{E}^{\mathbf{x}^*} \mathbf{x}_{\alpha}^{\delta}\|^2 + \mbox{tr }[\mathcal{C}_{\alpha}^{\delta}],
\end{equation}
where $\mathbf{x}_{\alpha}^{\delta}$ and $\mathcal{C}_{\alpha}^{\delta}$ are denoting the posterior mean and covariance, respectively. It suffices to bound the bias $b_{\mathbf{x}^*}(\alpha)= \|\mathbf{x}^* - \mathbb{E}^{\mathbf{x}^*}\mathbf{x}_{\alpha}^{\delta} \|$ and the posterior spread $\mbox{tr }[\mathcal{C}_{\alpha}^{\delta}]$, as the estimation variance $\mathbb{E}^{\mathbf{x}^*}\|\mathbf{x}_{\alpha}^{\delta} - \mathbb{E}^{\mathbf{x}^*} \mathbf{x}_{\alpha}^{\delta}\|^2$ is always bounded by the posterior spread.
To bound the posterior spread, it is shown that the necessary methods given in \cite{agapioumathe2018, mathe2018} can be extended to block matrix operators. First, the concept of index functions is discussed.
\begin{definition}
(index function) A function $\psi : (0, \infty) \to \mathbb{R}^+$ is called an index function if it is a continuous non-decreasing function with $\psi (0) = 0$.
\end{definition}
To be able to compare index functions the following partial ordering is defined.
\begin{definition}
(partial ordering for index functions) Given two index functions $g,h$ we write $g \prec h$ if the function $h(t)/g(t)$ is an index function (which means $h$ tends to zero faster than $g$).
\end{definition}
In addition to ordering index function, also self-adjoint operators need to be partially ordered.
\begin{definition}(partial ordering for self-adjoint operators) Let $\mathcal{G}$ and $\mathcal{G}'$ be bounded self-adjoint opertors in some Hilbert space $X \times Y$. We say that $\mathcal{G} \leq \mathcal{G}'$ if for all $\mathbf{x} \in X \times Y$ the inequality $\langle \mathcal{G} \mathbf{x}, \mathbf{x} \rangle \leq \langle \mathcal{G}' \mathbf{x}, \mathbf{x} \rangle$ holds.
\end{definition}
The last definition of interest states operator concavity for self-adjoint operators.
\begin{definition}
Let $f: [0,a] \to \mathbb{R}^+$ be a continuous function. It is called operator concave if we have for any pair $\mathcal{G}, \mathcal{H} \geq 0$ of self-adjoint operators with spectra in $[0,a]$
$$
f \left( \frac{\mathcal{G}+\mathcal{H}}{2} \right) \geq \frac{f(\mathcal{G}) +f(\mathcal{H})}{2}.
$$
\end{definition}
The theorem stated below is the range inclusion theorem, which also holds for block operator matrices.
\begin{theorem}(Douglas' Range Inclusion Theorem; Theorem 1 in \cite{Douglas})
Let the operators $\mathcal{S}, \mathcal{T}: X \times Y \to X \times Y$ be bounded and act between Hilbert spaces. Then the following statements are equivalent:
\begin{enumerate}
\item $\mathcal{R}(\mathcal{S}) \subset \mathcal{R}(\mathcal{T})$
\item $\mathcal{S}\mathcal{S}^* \leq C^2 \mathcal{T} \mathcal{T}^*$, for some $C \geq 0$
\item there exists an bounded operator $\mathcal{R}: X \times Y \to X \times Y$ with $||\mathcal{R}|| \leq C$, such that $\mathcal{S} = \mathcal{T}\mathcal{R}$.
\end{enumerate}
\end{theorem}
 
Before continuing we introduce the functions
\begin{equation}
\label{theta}
\Theta(t) = \Theta_{\psi}(t)  := \sqrt{t} \psi (t), \quad t > 0, \quad s_{\alpha}(t) := \frac{\alpha}{\alpha + t}, \quad \alpha >0,
\end{equation}
which will be needed in the further analysis. The main concept of the convergence analysis in \cite{mathe2018} is the fulfillment of two assumptions, namely a link and a source condition, which will be given next. With these two assumptions one is able to bound the squared bias and the posterior spread, and therefore, the squared posterior contraction.
\begin{assumption}(link condition)\label{ass:linkcond}
There is an index function $\psi$, and there are constants $0 < \underline{m} \leq \overline{m} < \infty$ such that
\begin{equation} \label{linkcond}
\underline{m}\| \psi (\mathcal{C}_0) \mathbf{x}\|_{W'\times Y} \leq \| \Sigma^{1/2} \mathcal{G} \mathbf{x}\|_{U\times X} \leq \overline{m} \| \psi ( \mathcal{C}_0) \mathbf{x} \|_{W'\times Y}, \quad \mathbf{x} \in U \times X. 
\end{equation}
In addition, with the function $\Theta$ from (\ref{theta}) the function
$$
f_0(s) := \left( (\Theta^2)^{-1} (s) \right)^{-1/2}, ~~ s >0,
$$
has an operator concave square $f_0^2$.
\end{assumption}
For the following analysis we define
\begin{equation}
\varphi_0(t) := \sqrt{t}, ~~ t > 0
\end{equation}
and state the following proposition, which will allow us to bound the bias.
\begin{proposition}\label{prop:normbound} (Proposition 1 in \cite{mathe2018})
Under Assumption~\ref{ass:linkcond}  we have that $\mathcal{R}(\mathcal{C}_0^{1/2}) = \mathcal{R}(f_0(\mathcal{H}))$. Mainly we have that the operator $f_0(\mathcal{H}) \varphi_0(\mathcal{C}_0)^{-1}$ is norm bounded by $\overline{m}$.
\end{proposition}

As mentioned before, the second assumption of importance is a source condition, which will be stated next.
\begin{assumption}(source set)\label{ass:sourceset}
There is an index function $\varphi$ such that
\begin{equation}
\mathbf{x}^* \in \mathcal{S}_{\varphi} := \{\mathbf{x}, \quad \mathbf{x} = \varphi(\mathcal{C}_0)\mathbf{v}, ~||\mathbf{v}|| \leq 1 \}.
\end{equation}
\end{assumption}
By using Proposition~\ref{prop:normbound} we can bound the bias as
\begin{eqnarray}
b_{\mathbf{x}^*}(\alpha) &\leq \frac{1}{\underline{m}} ||f_0(\mathcal{H}) s_{\alpha}(\mathcal{H}) \varphi_0(\mathcal{C}_0)^{-1}\varphi(\mathcal{C}_0)|| \nonumber\\
&= \frac{1}{\underline{m}} ||s_{\alpha}(\mathcal{H}) f_0(\mathcal{H}) \varphi_0(\mathcal{C}_0)^{-1}\varphi(\mathcal{C}_0)||. \label{estbias0}
\end{eqnarray}
Now, with the assumptions made, one can give an upper bound for the squared bias.
\begin{proposition}\label{prop:bias}(Proposition 3 in \cite{mathe2018})
Suppose that either $\varphi_0 \prec \varphi \prec \theta$, and the function 
\begin{equation*}
g^2 (t) := \left( \frac{\varphi}{\varphi_0} \right)^2 \left( \left( \theta^2 \right)^{-1}(t) \right), \quad t >0,
\end{equation*}
is operator concave, or $1  \prec \varphi \prec \varphi_0$ and $\varphi$ is operator concave. Under Assumptions~\ref{ass:linkcond} and \ref{ass:sourceset} we have that
\begin{equation}\label{est_bias}
b_{\mathbf{x}^*}(\alpha) \leq \frac{\overline{m}}{\underline{m}}||s_{\alpha}(\mathcal{H}) \varphi (f_0^2(\mathcal{H}))||.
\end{equation}
\end{proposition}
With the help of Assumption~\ref{ass:linkcond}, one can bound the posterior spread as stated in the next proposition.
\begin{proposition}\label{prop:postspread} (Proposition 6 in \cite{mathe2018}) Under Assumption~\ref{ass:linkcond} we have that 
\begin{equation}
\mbox{tr}[\mathcal{C}_{\alpha}^{\delta}] \leq \frac{\delta^2}{\underline{m}^2} \mbox{tr}[(\alpha + \mathcal{H})^{-1}f_0^2(\mathcal{H})].
\end{equation}
\end{proposition}

Under certain circumstances, see \cite{mathe2018} and the references therein, we can rewrite the upper bound for the bias as
\begin{equation}
b_{\mathbf{x}^*} (\alpha) \leq \frac{\overline{m}}{\underline{m}} \varphi(f_0^2(\alpha)), \quad \alpha > 0,
\end{equation}
using the fact that $s_{\alpha}(t)t \leq \alpha$, for $t,\alpha >0$.

Combination of the results from Proposition~\ref{prop:bias} and \ref{prop:postspread} leads to a bound on the squared posterior contraction as stated in the following theorem.
\begin{theorem}(Theorem 3 in \cite{mathe2018})
Supposed that Assumption~\ref{ass:linkcond} and Assumption~\ref{ass:sourceset} hold for some index functions $\psi$ and $\varphi$ and additionally, the assumptions in Proposition~\ref{prop:bias} and holds true, then
\begin{equation}\label{SPC_full}
\mbox{SPC}(\alpha, \delta) \leq \frac{1}{\underline{m}^2} 
\Bigl(\overline{m}^2 ||s_{\alpha}(\mathcal{H}) \varphi (f_0^2(\mathcal{H}))|| 
+ \delta^2 \mbox{ tr }[(\alpha + \mathcal{H})^{-1} f_0^2(\mathcal{H})]\Bigr), \quad \alpha, \delta >0.
\end{equation}
\end{theorem}

\begin{remark}\label{rem:lb}
If only the lower bound in \eqref{linkcond} holds, with  \eqref{estbias0} and under a source condition one still obtains the estimate
\begin{equation}\label{SPC_half}
\mbox{SPC}(\alpha, \delta) = \frac{1}{\underline{m}^2} 
\Bigl(||s_{\alpha}(\mathcal{H}) f_0(\mathcal{H})\varphi_0(\mathcal{C}_0)^{-1}\varphi(\mathcal{C}_0)|| 
+ \delta^2 \mbox{ tr }[(\alpha + \mathcal{H})^{-1} f_0^2(\mathcal{H})]\Bigr),
\end{equation}
for $\alpha, \delta >0.$
\end{remark}

\subsection{Fulfillment of the link condition for the all-at-once-formulation.}
A straightforward example that always allows to fulfill the link condition, is as follows. Choosing the prior covariance to be given by
\begin{equation}
\label{choicecov}
\mathcal{C}_0 := \mathcal{G}^* \mathcal{G},
\end{equation}
and the link condition is always fulfilled by choosing the link function to be given by
\begin{equation}
\psi (t) := t^{1/2}.
\end{equation}
Assume $\mathcal{G}:X \times Y \to S \times T$, then $\mathcal{G}^* : S \times T \to X \times Y$ and $\mathcal{G}^*\mathcal{G} : X \times Y \to X \times Y$ which yields
\begin{equation}
\| \psi (\mathcal{G}^* \mathcal{G}) \mathbf{x}\|_{X \times Y} = \| \mathcal{G}\mathbf{x} \|_{S \times T}, \quad \mbox{ for all } \mathbf{x} \in X \times Y,
\end{equation}
and the link condition is fulfilled. However, since $\mathcal{G}^*\mathcal{G}$ is ill-posed, \eqref{choicecov} will not be a good regularizer, so we look at further possible choices for our two examples in the following Section \ref{sec:priors}.

In particular, for the backwards heat equation, the severe ill-posedness is provibitive for \eqref{choicecov}.
In the Section \ref{sec:heuristicprior} below, we therefore derive a prior that is motivated by the link condition but mach easier to handle than \eqref{choicecov}.

\section{Choice of joint priors.} \label{sec:priors}

\subsection{Block diagonal priors satisfying unilateral link estimates.}
In this section, we investigate to which extent the link condition \eqref{linkcond} can be satisfied by block diagonal operators $\psi(\mathcal{C}_0)$, as these are convenient for applying $\psi^{-1}$ in order to obtain the prior covariance $\mathcal{C}_0$ itself again in block diagonal form, which in its turn is useful for numerical computations. For this purpose we again focus on our two prototypical examples from Sections \ref{sec:is} and \ref{sec:bh}. We will consider the upper and lower bound in \eqref{linkcond} separately and, assuming that $\Sigma$ is an isomorphism on $Y$, without loss of generality, set it to the identity, i.e., we will investigate fulfillment of 
\begin{equation} \label{linkcond_lb}
\| \psi (\mathcal{C}_0) \mathbf{x}\|_{U\times X} \leq \frac{1}{\underline{m}} \| \mathcal{G} \mathbf{x}\|_{W'\times Y} \mbox{ for all } \mathbf{x} = (u,\theta)\in U\times X,
\end{equation}
\begin{equation} \label{linkcond_ub}
\| \mathcal{G} \mathbf{x}\|_{W'\times Y} \leq \overline{m} \| \psi ( \mathcal{C}_0) \mathbf{x} \|_{U\times X} \mbox{ for all } \mathbf{x} = (u,\theta)\in U\times X,
\end{equation}
for $\psi (\mathcal{C}_0)$ of the form 
\begin{equation}\label{diag_prior}
\psi(\mathcal{C}_0) = \left( \begin{array}{cc} \mathcal{B} & 0 \\ 0 & \mathcal{D} \end{array} \right).
\end{equation}

For the inverse source problem from Section \ref{sec:is} we have 
\[
\|\mathcal{G}(u,\theta)\|_{W'\times Y} = \|\mathcal{A}u-\theta\|_{L^2}+\|u\|_{L^2}
\]
and
\[ 
\| \psi ( \mathcal{C}_0) (u,\theta) \|_{U\times X} = 
\|\mathcal{A}\mathcal{B}u\|_{L^2}+\|\mathcal{D}\theta\|_{L^2}
\]
so setting $\theta$ and $u$ to zero separately we get the following two necessary conditions 
\begin{eqnarray}\label{B_is_ub}
&& \|\mathcal{A}u\|_{L^2}\leq \overline{m}\|\mathcal{A}\mathcal{B}u\|_{L^2} \mbox{ for all } u\in U= H_0^1(\Omega) \cap H^2(\Omega)\,, \nonumber\\
&& \mbox{ i.e., $\mathcal{B}:U \to U$ bijective with bounded inverse $\mathcal{B}^{-1}$}
\end{eqnarray}
and
\begin{eqnarray}\label{D_is_ub}
&& \|\theta\|_{L^2}\leq \overline{m}\|\mathcal{D}\theta\|_{L^2} \mbox{ for all } \theta\in X= L^2(\Omega)\,, \nonumber\\
&& \mbox{ i.e., $\mathcal{D}:L^2(\Omega)\to L^2(\Omega)$ bijective with bounded inverse $\mathcal{D}^{-1}$}
\end{eqnarray}
for \eqref{linkcond_ub}. On the other hand, an easy application of the triangle inequality together with boundedness of $\mathcal{A}^{-1}:L^2(\Omega)\to L^2(\Omega)$ implies sufficiency of \eqref{B_is_ub}, \eqref{D_is_ub} for \eqref{linkcond_ub}.
\\
The lower bound \eqref{linkcond_lb}, with $\theta=\mathcal{A}u$ implies 
\[
\|\mathcal{A}\mathcal{B}u\|_{L^2}+\|\mathcal{D}\mathcal{A}u\|_{L^2} \leq \frac{1}{\underline{m}}\|u\|_{L^2} \mbox{ for all }u\in U= H_0^1(\Omega) \cap H^2(\Omega)\,,
\]
thus since $\|\mathcal{D}\mathcal{A}\|_{L^2\to L^2} = \|\mathcal{A}\mathcal{D}^*\|_{L^2\to L^2}$
we get the necessary conditions
\begin{equation}\label{BD_is_lb}
\mbox{$\mathcal{B}:L^2(\Omega)\to H_0^1(\Omega) \cap H^2(\Omega)$ and $\mathcal{D}^*:L^2(\Omega)\to H_0^1(\Omega) \cap H^2(\Omega)$ bounded} 
\end{equation}
for \eqref{linkcond_lb}. To see sufficiency of \eqref{BD_is_lb} together with boundedness of $\mathcal{D}:L^2(\Omega)\to L^2(\Omega)$ for \eqref{linkcond_lb}, we set $f=\mathcal{A}u-\theta$ to obtain
\begin{eqnarray*}
\|\psi ( \mathcal{C}_0) (u,\theta) \|_{U\times X} &=& 
\|\mathcal{A}\mathcal{B}u\|_{L^2}+\|\mathcal{D}(\mathcal{A}u-f)\|_{L^2}\\
&\leq&(\|\mathcal{B}\|_{L^2\to H^2}+\|\mathcal{D}^*\|_{L^2\to H^2})\|u\|_{L^2}
+\|\mathcal{D}\|_{L^2\to L^2})\|f\|_{L^2}\\
&\leq& \frac{1}{\underline{m}}(\|f\|_{L^2}+\|u\|_{L^2}) 
= \frac{1}{\underline{m}} \|\mathcal{G}(u,\theta)\|_{W'\times Y}
\end{eqnarray*}
with $\underline{m}:=\frac{1}{\max\{\|\mathcal{B}\|_{L^2\to H^2}+\|\mathcal{D}^*\|_{L^2\to H^2},\,
\|\mathcal{D}\|_{L^2\to L^2}\}}$.
\\
Unfortunately, conditions \eqref{B_is_ub} and \eqref{BD_is_lb} are contradictory, so the full squared posterior contraction estimate \eqref{SPC_full} cannot be applied.
Still, with, e.g. $\psi(\mathcal{C}_0) = \left( \begin{array}{cc} \mathcal{A}^{-s} & 0 \\ 0 & \mathcal{A}^{-p} \end{array} \right)$ and $s,p\geq1$, which satisfies \eqref{linkcond_lb}, we obtain \eqref{SPC_half}.
\\
For the backwards heat problem from Section \ref{sec:is} with 
\[
\|\mathcal{G}(u,\theta)\|_{W'\times Y} = \|\mathcal{A}^{-1/2}(\partial_t+\mathcal{A})u+\mathcal{A}^{1/2}\theta\|_{L^2(L^2)}+\|u(T)+\theta\|_{L^2}
\]
and
\[ 
\| \psi ( \mathcal{C}_0) (u,\theta) \|_{U\times X} = 
\|\mathcal{A}^{-1/2}(\partial_t+\mathcal{A})\mathcal{B}u\|_{L^2(L^2)}+\|\mathcal{A}^{1/2}\mathcal{D}\theta\|_{L^2}
\]
we see, by considering the two special cases $\theta=0$ and $u=0$, that a necessary condition for \eqref{linkcond_ub} is 
\begin{equation}\label{BD_bh_ub}
\mbox{$\mathcal{B}^{-1}:U_0\to U_0$ and $\mathcal{D}^{-1}:H_0^1(\Omega)\to H_0^1(\Omega) $ bounded.}
\end{equation}
Indeed, as can be seen by the triangle inequality and fact that that $\|u(T)\|_{L^2}\leq \|\mathcal{A}^{-1/2}(\partial_t+\mathcal{A})u\|_{L^2(L^2)}$ (cf. \eqref{innprodU0}) condition \eqref{BD_bh_ub} is also sufficient for \eqref{linkcond_ub}.
\\
The lower bound \eqref{linkcond_lb} is again more challenging to obtain. With the particular choices $\theta=0$ on one hand and $(\partial_t+\mathcal{A})u=-\mathcal{A}\theta$, i.e., $u(t)=-(I-e^{-t\mathcal{A}})\theta$ implying $u(T)+\theta = e^{-T\mathcal{A}}\theta$ on the other hand, it yields the necesssary conditions 
\begin{eqnarray*}
&&\|\mathcal{B}u\|_{U_0} \leq \frac{1}{\underline{m}} \Bigl(\|u\|_{U_0}+\|u(T)\|_{L^2}\Bigr)\leq
\frac{2}{\underline{m}} \|u\|_{U_0}  \mbox{ for all }u\in U_0,\\
&&\mbox{i.e., $\mathcal{B}:U_0\to U_0$ bounded}\nonumber
\end{eqnarray*}
and 
\begin{eqnarray*}
&&\|\mathcal{B}(\ulI-e^{-\cdot\mathcal{A}})\theta\|_{U_0}+
\|\mathcal{D}\theta\|_{H_0^1} \leq \frac{1}{\underline{m}} \|e^{-T\mathcal{A}}\theta\|_{L^2} \mbox{ for all }\theta\in X,\\
&&\mbox{i.e., $\mathcal{D} e^{T\mathcal{A}}:L^2(\Omega)\to H_0^1(\Omega)$ bounded.} \nonumber
\end{eqnarray*}
To obtain sufficient conditions for \eqref{linkcond_lb}, with $f=(\partial_t+\mathcal{A})u+\mathcal{A}\theta$, $g=u(T)+\theta$, so that 
$\|\mathcal{G}(u,\theta)\|_{W'\times Y} = \|f\|_{L^2(H^{-1})}+\|g\|_{L^2}$, 
$u(t)=\int_0^t e^{-(t-s)\mathcal{A}} f(s)\, ds -(I-e^{-t\mathcal{A}})\theta$, hence we can express 
$u$ and $\theta$ via $f$ and $g$ as 
\begin{eqnarray*}
u(t)&=&\int_0^t e^{-(t-s)\mathcal{A}} f(s)\, ds -(I-e^{-t\mathcal{A}}) (e^{T\mathcal{A}} g-\int_0^T e^{s\mathcal{A}} f(s)\, ds),\\
\theta&=& e^{T\mathcal{A}} g-\int_0^T e^{s\mathcal{A}} f(s)\, ds ,
\end{eqnarray*}
and, therefore ,
\begin{align*}
 \| \psi ( \mathcal{C}_0) (u,\theta) \|_{U\times X} &=
\|\mathcal{B}u\|_{U_0}+\|\mathcal{D}\theta\|_{H_0^1}\\
& \leq \|\mathcal{B}\|_{L^\infty(L^2)\to U_0}
\|\int_0^\cdot e^{-(\cdot-s)\mathcal{A}} f(s)\, ds\|_{L^\infty(L^2)}\\
& \quad+ \Bigl(\|\mathcal{B}e^{T\mathcal{A}}\|_{L^\infty(L^2)\to U_0}
\|\ulI-e^{-\cdot\mathcal{A}}\|_{L^2\to L^\infty(L^2)} 
+\|\mathcal{D}e^{T\mathcal{A}}\|_{L^2\to H_0^1}\Bigr)\\
& \qquad\qquad\cdot\Bigl(\|g\|_{L^2}+ \|\int_0^T e^{-(T-s)\mathcal{A}} f(s)\, ds\|_{L^2}\Bigr)\,,
\end{align*}
where $\|\ulI-e^{-\cdot\mathcal{A}}\|_{L^2\to L^\infty(L^2)}\leq1$
and by self-adjoinedness of $\mathcal{A}$ as well as the Cauchy-Scharz inequality 
\begin{align*}
& \|\int_0^\cdot e^{-(\cdot-s)\mathcal{A}} f(s)\, ds\|_{L^\infty(L^2)} \\
&=\sup_{t\in[0,T]\, v\in L^2\setminus\{0\}} \tfrac{1}{\|v\|_{L^2}}
\int_0^t \int_\Omega e^{-(t-s)\mathcal{A}} \mathcal{A}^{1/2} v\, \mathcal{A}^{-1/2} f(s)\, dx\, ds \\
& \leq \sup_{t\in[0,T]\, v\in L^2\setminus\{0\}} \tfrac{1}{\|v\|_{L^2}}
\left(\int_0^t \sum_{n=1}^\infty  e^{-2(t-s)\frac{1}{\mu_n}} \frac{1}{\mu_n} \langle v,\phi_n\rangle_{L^2}^2 \,ds\right)^{1/2}\ \|\mathcal{A}^{-1/2}f\|_{L^2(L^2)} \\
&\leq\frac12 \|f\|_{L^2(H^{-1})}
\end{align*}
with an eigensystem $(\mu_n,\phi_n)_{n\in\mathbb{N}}\subseteq\mathbb{R}\times L^2(\Omega)$ of $\calA^{-1}$; likewise \\
$\|\int_0^T e^{-(T-s)\mathcal{A}} f(s)\, ds\|_{L^2}\leq\frac12 \|f\|_{L^2(H^{-1})}$.
Since these estimates are sharp when having to hold for arbitrary $f\in L^2(H^{-1})$, we can conclude that validity of \eqref{linkcond_lb} is equivalent to 
\begin{equation}\label{BD_bh_lb}
 \mbox{$\mathcal{B}e^{T\mathcal{A}}:L^\infty(L^2)\to U_0$ and $\mathcal{D} e^{T\mathcal{A}}:L^2(\Omega)\to H_0^1(\Omega)$ bounded.} 
\end{equation}
Again the conditions for \eqref{linkcond_lb} and \eqref{linkcond_ub} are unfortunately contradictory. Still, for example the choice 
$\psi(\mathcal{C}_0) = \left( \begin{array}{cc} \mathcal{A}^{1/2}e^{-T\mathcal{A}}\int_0^\cdot e^{-(\cdot-s)\mathcal{A}} \cdot(s)\, ds & 0 \\ 0 & \mathcal{A}^{-1/2}e^{-T\mathcal{A}}\end{array} \right)$ satisfies \eqref{linkcond_lb} and therewith \eqref{SPC_half}.

\subsection{Heuristic choice of $\mathcal{C}_0$ for the backwards heat problem.}\label{sec:heuristicprior}
To find out more about suitable operators fulfilling the link condition besides the trivial choice one attempt is made in the setting of the backwards heat problem. To do so, the norms are compared and so, the operator $\psi(\mathcal{C}_0)$ can be given by
\[
\psi(\mathcal{C}_0) = \left( \begin{array}{cc} I & \underline{I} - e^{-\mathcal{A}t} \\ \mathcal{A}^{-1} \delta_T + B & \mathcal{A}^{-1/2} \end{array} \right),
\]
where $\mathcal{A} := - \Delta$, $Bu := \int_0^T u(t) \mbox{ dt}$ and $(\underline{I}\theta)(t) = \theta$.
This is a symmetric operator from $U_0 \times H_0^1(\Omega) \to U_0 \times H_0^1(\Omega)$ as
\begin{align*}
 ((I - e^{-\mathcal{A}t})u, v)_{U_0} &= \int_0^T \int_{\Omega}  \mathcal{A}^{-1/2} (\partial_t + \mathcal{A})(\underline{I}-e^{-\mathcal{A}t})u \mathcal{A}^{-1/2}(\partial_t + \mathcal{A})v \mbox{ dx} \mbox{ dt}  \\
&=  \int_0^T \int_{\Omega} \mathcal{A}^{1/2}u \mathcal{A}^{-1/2}(\partial_t + \mathcal{A})v \mbox{ dx} \mbox{ dt}\\
&= \int_0^T \int_{\Omega}\mathcal{A}^{1/2}u \mathcal{A}^{-1/2}\partial_t v \mbox{ dx} \mbox{ dt} + \int_0^T\int_{\Omega} \mathcal{A}^{1/2}u \mathcal{A}^{1/2}v \mbox{ dx} \mbox{ dt} \\
 &=  \int_{\Omega} \mathcal{A}^{1/2} u \mathcal{A}^{-1/2} v(T) \mbox{ dx} + \int_{\Omega} \mathcal{A}^{1/2}u \mathcal{A}^{1/2} \int_0^T v(t) \mbox{ dt} \mbox{ dx} ,
 \end{align*}
 \begin{align*}
 (u, (\mathcal{A}^{-1} \delta_T + B)v)_{H_0^1(\Omega)} &= \int_{\Omega} \mathcal{A}^{1/2}u \mathcal{A}^{1/2}(\mathcal{A}^{-1} \delta_T + B)v \mbox{ dx} \\
&= \int_{\Omega} \mathcal{A}^{1/2} u \mathcal{A}^{-1/2}v(T) \mbox{ dx} + \int_{\Omega} \mathcal{A}^{1/2} u \mathcal{A}^{1/2}\int_0^Tv(t) \mbox{ dt}.
\end{align*}
Now the norm of $\mathcal{G}$ and $\psi(\mathcal{C})$ are computed with the help of a rewriting of $\psi(\mathcal{C})$ as
\begin{eqnarray}
\psi(\mathcal{C}_0) &=& \left( \begin{array}{cc} I & \underline{I} - e^{-\mathcal{A}t} \\ \mathcal{A}^{-1} \delta_T + B & \mathcal{A}^{-1/2} \end{array} \right) \label{Ctil}\\
&=& \underbrace{\left( \begin{array}{cc} I & 0 \\ 0 & \mathcal{A}^{-1/2} \end{array} \right)}_{=: \tilde{\mathcal{A}}} \underbrace{\left( \begin{array}{cc} I & \underline{I}-e^{- \mathcal{A}t} \\ \mathcal{A}^{-1/2} \delta_T & \mathcal{A}^{-1/2} \end{array} \right)}_{\tilde{\mathcal{C}}} + \underbrace{\left( \begin{array}{cc} 0 & 0 \\ B & \mathcal{A}^{-1/2} - \mathcal{A}^{-1} \end{array} \right)}_{=: \mathcal{R}}.
\nonumber
\end{eqnarray}
Instead of $\psi(\mathcal{C}_0)$ now we use $\tilde{\mathcal{C}}$ under the norm. The norms compute as follows for $x = (u, \theta) \in U_0 \times H_0^1$
\begin{eqnarray*}
 \| \mathcal{G} \textbf{x} \|_{L^2(H^{-1}) \times L^2}^2 &=& \| (\partial_t +\mathcal{A})u + \mathcal{A}\theta||_{L^2(H^{-1})}^2 + \| \delta_T u + I\theta \|_{L^2}^2  \\
 &=&\| (\partial_t + \mathcal{A})u\|_{L^2(H^{-1})}^2 + \| \mathcal{A}\theta \|_{L^2(H^{-1})}^2 + 2((\partial_t + \mathcal{A})u, \mathcal{A}\theta)_{L^2(H^{-1})} \\
 &&+ \| \delta_T u\|_{L^2}^2 + \| I \theta \|_{L^2}^2 + 2(\delta_T u, I\theta)_{L^2}  \\
 &=& \int_0^T \int_{\Omega} |\mathcal{A}^{-1/2} (\partial_t + \mathcal{A})u|^2 \mbox{ dx} \mbox{ dt} + \int_0^T \int_{\Omega}|\mathcal{A}^{1/2} \theta |^2 \mbox{ dx}\mbox{ dt} \\
 &&+  2 \int_0^T \int_{\Omega} \mathcal{A}^{-1/2} (\partial_t + \mathcal{A})u \mathcal{A}^{1/2} \theta \mbox{ dx} \mbox{ dt}+\int_{\Omega} |u(T)|^2 \mbox{ dx} \\
 &&+ \int_{\Omega} |\theta|^2 \mbox{ dx} + 2\int_{\Omega} u(T) I\theta \mbox{ dx},
\end{eqnarray*}

\begin{eqnarray*}
 \| \tilde{\mathcal{C}} \textbf{x}\|_{U_0 \times H_0^1}^2 &=& \| u + (\underline{I}-e^{-\mathcal{A} t} \theta )||_{U_0}^2 + \| \mathcal{A}^{-1/2}\delta_T u + \mathcal{A}^{-1/2}\theta \|_{H_0^1} ^2\\
 &=&\| u\|_{U_0}^2 + \| (I-e^{-\mathcal{A} t})\theta \|_{U_0} ^2+ 2(u, (\underline{I}- e^{- \mathcal{A} t})\theta)_{U_0} \\
 &&+ \| \mathcal{A}^{-1/2}\delta_T u\|_{H_0^1}^2 + \| \mathcal{A}^{-1/2} \theta \|_{H_0^1}^2 + 2(\mathcal{A}^{-1/2} \delta_T u, \mathcal{A}^{-1/2} \theta)_{H_0^1} \\
 &=& \int_0^T\int_{\Omega} |\mathcal{A}^{-1/2} (\partial_t + \mathcal{A})u|^2 \mbox{ dx} \mbox{ dt} \\
 &&+ \int_0^T \int_{\Omega} |\mathcal{A}^{-1/2} (\partial_t + \mathcal{A}) (I - e^{-\mathcal{A} t})\theta |^2 \mbox{ dx} \mbox{ dt} \\
 &&+ 2\int_0^T \int_{\Omega} \mathcal{A}^{-1/2} (\partial_t + \mathcal{A}) u \mathcal{A}^{-1/2} (\partial_t + \mathcal{A}) (\underline{I} - e^{-\mathcal{A} t})\theta \mbox{ dx} \mbox{ dt} \\
 &&+ \int_{\Omega} | \mathcal{A}^{1/2} \mathcal{A}^{-1/2} \delta_T u|^2 \mbox{ dx} \\
 &&+ \int_{\Omega} |\mathcal{A}^{1/2} \mathcal{A}^{-1/2} \theta |^2 \mbox{ dx} + 2 \int_{\Omega} \mathcal{A}^{1/2} \mathcal{A}^{-1/2} \delta_T u \mathcal{A}^{1/2}\mathcal{A}^{-1/2} \theta \mbox{ dx} 
\\
 &=& \int_0^T\int_{\Omega} |\mathcal{A}^{-1/2} (\partial_t + \mathcal{A})u|^2 \mbox{ dx} \mbox{ dt} + \int_0^T \int_{\Omega} |\mathcal{A}^{1/2} \theta |^2 \mbox{ dx} \mbox{ dt} \\
 &&+ 2\int_0^T \int_{\Omega} \mathcal{A}^{-1/2} (\partial_t + \mathcal{A}) u \mathcal{A}^{1/2} \theta \mbox{ dx} \mbox{ dt} + \int_{\Omega}   |u(T)|^2 \mbox{ dx} \\
 &&+ \int_{\Omega}  |\theta |^2 \mbox{ dx} + 2 \int_{\Omega} v(T) u \mbox{ dx}.
\end{eqnarray*}
This shows that $\| \tilde{\mathcal{C}} \mathbf{x} \|^2 = \|\mathcal{G} \mathbf{x} \|^2$. 
Motivated by this, we use $\tilde{ \mathcal{C}}$ in place of $\mathcal{C}$ as a prior covariance matrix. Note that estimating the remainder $\|\mathcal{R} \mathbf{x}\|$ by $\|\mathcal{G} \mathbf{x}\|$ unfortunately does not seem to be possible, since in $\mathcal{G} \mathbf{x}$ terms containing $u$ and terms containing $\theta$ may possibly cancel. Nevertheless, since $\mathcal{R}$ contains negative order differential operators as compared to $\tilde{ \mathcal{C}}$, we can regard it as a perturbation of the latter and therefore skip it for computational purposes.

In Subsection \ref{sec:heuristic} below we demonstrate numerically that block diagonal of $\tilde{\mathcal{C}}$ as a prior yields reasonable results in the reconstruction.

\subsection{Priors for the inverse source problem.}
The choice of joint priors in this work relies on priors which are already well known in the literature and also implemented in hippylib, see, e.g. \cite{stuart:dashti}\cite{VillaPetraGhattas2018}. We restrict ourselves to normal priors with the covariance operator given by
\begin{equation}
C_p= (- \gamma \Delta + \delta I)^{-n}
\end{equation}
where $n$ most often take the values $1$ or $2$. For joint priors for the parameter and the state, priors of the same type have shown to be useful in the numerical experiments in this study. 

\subsection{Prior for the state variable of the backwards heat problem.} \label{sec:semigroup-prior}
For the backwards heat problem we require a prior which is suitable for the initial condition, but in addition, we also need to find a prior for the state variable for every $t \in [0,T]$. To do so, the heat equation is analyzed. 
At first the solution to the heat equation is considered and the homogeneous case $f=0$ is examined. 
In case of bounded $\Omega$ with, e.g., homogeneous Dirichlet boundary conditions and denoting by $\calA=-\Delta$ the Laplace operator equipped with these boundary conditions, we can use semigroup theory to express $u$.
To use the semigroup theory some preparatory work has to be done. 
\begin{definition}(Section 7.4.1 in \cite{evans:PDE})
Let $X$ be a real Banach space. A family of linear, bounded operators $\{S(t)\}_{t \geq 0}$, $S(t):X \to X$  is called a semigroup iff
\begin{equation*}
S(0) = I \quad \mbox{and} \quad S(s+t) = S(s)S(t) \quad \forall t,s \geq 0.
\end{equation*}
\end{definition}
As already assumed the solution of the heat equation with initial data $\theta$ can be written in terms of semigroups using the so called heat semigroup
\begin{equation}
S(t)\theta:=u(t)   \mbox{ where } \left\{\begin{array}{rcl}u_t&=&\Delta u \mbox{ in }\Omega\times(0,\infty)\\ u&=&0\mbox{ on } \partial\Omega\times(0,\infty)\\ \quad u&=&\theta  \mbox{ on }\Omega\times\{0\}\end{array}\right.
\end{equation}
with $X= L^2(\Omega)$.
Clearly, $S(t)$ fulfills the conditions to be a semigroup. To investigate this equation the definition of the infinitesimal generator is needed.
\begin{definition}(Section 7.4.1 in \cite{evans:PDE})
Let $\{S(t)\}_{t\geq 0}$ be a semigroup on a Banach space $X$. Write 
\begin{equation*}
D(A) := \left\{ u \in X: \lim_{t \to 0^+} \frac{S(t)u -u}{t} \mbox{ exists in } X \right\} 
\end{equation*}
and
\begin{equation*}
Au := \lim_{t \to 0+} \frac{S(t)u -u}{t} \quad (u \in D(A)).
\end{equation*}
We call $A: D(A) \to X$ the infinitesimal generator of the semigroup $\{S(t)\}_{t \geq 0}$, $D(A)$ is the domain of $A$.
\end{definition}
It can be shown, that the generator of the heat semigroup is defined by the Laplace operator with homogeneous Dirichet boundary conditions $\calA$. 
Therewith, we can write $u$ as 
\begin{equation} \label{u_semigr}
u(t,x) = (S(t)\theta)(x) = (e^{-t\calA} \theta)(x).
\end{equation}

In terms of finding a prior for the state $u(t,x)$ at any time $t$, this is a convenient setting as we can make use of the prior measure for $\theta(x)$. We denote the measure for the parameter $\theta$ by $\mu_0$ and given by
\begin{equation}
\mu_0 = \mathcal{N}(m_0, C_p) 
\end{equation}
Using the properties of the normal distribution one can calculate the resulting prior for $u(t,x)$. In addition the fact that $e^{-t\calA}$ is a self-adjoint operator is used. Then the prior for $u(t,x)$ is given by
\begin{equation}
\mu_{s} = \mathcal{N}( e^{-t\calA}m_0, e^{-t\calA} C_p e^{-t\calA}).
\end{equation}

In case $f\neq 0$, one can show that the solution to the heat equation changes to
\begin{equation}
u(t,x) = S(t)\theta(x) + \int_0^t S(t-s)f(s,x) \mbox{ ds} = e^{-t\calA} \theta(x) + \int_0^t e^{-(t-s)\calA}f(s,x) \mbox{ ds}.
\end{equation}
In the corresponding prior measure only the expectation value has changed, as the covariance operator is invariant under translation. This leads to the new prior
\begin{equation}
\mu_{s} = \mathcal{N}( e^{-t\calA}m_0 + \int_0^t e^{-(t-s)\calA}f(s,x)\mbox{ ds}, e^{-t\calA} C_p e^{-t\calA}).
\end{equation}

\section{Numerical Experiments.} \label{sec:num}
In this section the prototypical inverse problems discussed in the previous section will be solved numerically. To do so, first the theoretical considerations for solving the problems are described, while later the implementation is discussed. In this work, we used the software FEniCS \cite{finite:element} with the extension of hIPPYlib \cite{VillaPetraGhattas2018} which we extended further for the needs of our all-at-once formulations.
In this work the method for discretizing and solving the inverse problem is based on the methods described in \cite{thanh:inversebayes}. However, we change the computation of the Hessian including all the prior information we have in the all-at-once setting. 
The equation describing the relation between the data is given via
\begin{equation}
\obs u = y + \eta,
\end{equation}
where $\obs: U \to Y$ is the operator which maps the state variable into the observation space $Y$ and $\eta$ is the random variable describing the noise in this equation. We assume that $\eta$ is random white noise.

\subsection{Lagrangian method for computing the adjoint based Hessian and gradient.}
In this section we compute the adjoint based Hessian and gradient for a cost functional formulated in the all-at-once setting. This functional is then analyzed twice, first for the parameter and then also for the state variable. 

For the computation of the Hessian with respect to the parameter we consider the following minimization problem
\begin{align}
\min_{\theta} &J_\alpha(\theta)=\| \obs u-y \|_Y^2 +\frac{\alpha}{2} \mathcal{R}(u,\theta)\\
s.t. ~~&\mod(u,\theta) = 0\label{modutheta0}
\end{align}
with the regularization given by the prior distribution and written as
\begin{equation}
\mathcal{R}(u,\theta) = \langle u, C_1u \rangle_U + 2\langle u, C_2 \theta \rangle_U + \langle \theta, C_3 \theta \rangle_X ,
\end{equation}
with selfadjoint operators $C_1:U\to U$, $C_3:X\to X$.
Together, this leads to the following Lagrange functional of $u\in U$, $\theta\in X$, $p\in W$, 
\begin{align}
 \mathcal{L}_\alpha(u, \theta, p) :=& \frac{1}{2} \| \obs u-y \|_Y^2 + \frac{\alpha}{2}\big(\langle u, C_1u \rangle_U + 2\langle u, C_2 \theta \rangle_U + \langle \theta, C_3 \theta \rangle_X \big) \\
 &+\mod(u,\theta)p,
\end{align}
which is used in order to derive adjoint based methods for gradient and Hessian computation. 
We start with computing the first derivatives and for this purpose use the Riesz isomorphisms $I_R^U:U\to U'$, $u\mapsto (v\mapsto \langle u,v\rangle_U)$, $I_R^X:X\to X'$
\begin{align}
 \frac{\partial \mathcal{L}_\alpha}{\partial p} (u, \theta, p) &= \mod(u,\theta) \in W'\\
 \frac{\partial \mathcal{L}_\alpha}{\partial u} (u, \theta, p) &= 
I_R^U\Bigl(\obs^*(\obs u-y) + \alpha \left( C_1 u + C_2 \theta \right)\Bigl)  
+ \frac{\partial \mod}{\partial u}(u,\theta)^\star p \in U'\\
 \frac{\partial \mathcal{L}_\alpha}{\partial \theta} (u, \theta, p) &= 
\alpha I_R^X\left( C_2^* u + C_3 \theta \right) 
+ \frac{\partial \mod}{\partial \theta}(u,\theta)^\star p \in X'\,, 
\end{align}
where the derivative with respect to $p$ set to zero gives the forward problem. 
The following operators are now used as abbreviations of the derivatives of $\mod(u, \theta)$
\begin{equation}
-\frac{\partial \mod}{\partial u} (u, \theta) =: K:U\to W', \qquad 
\frac{\partial \mod}{\partial \theta} (u, \theta) =: L:X\to W',
\end{equation} 
along with their Banach space adjoints $K^\star:W\to U'$, $L^\star:W\to X'$ and their liftings to $U$ and $X$, respectively $K^*:=(I_R^U)^{-1}K^\star:W\to U$, $L^*:=(I_R^X)^{-1}L^\star:W\to X$.
Since we deal with linear models here, $K$ and $L$ are in fact independent of $u$ and $\theta$ and we can write $\mod(u,\theta)=-Ku+L\theta+f$.
Then, the derivative with respect to $u$ when set to zero leads to the adjoint problem with the adjoint linearization of the model defining the differential operator and the residual plus some additional term from the prior define the right hand side
\begin{equation}\label{adjointeq}
K^*p = \obs^*(\obs u - y) + \alpha (C_1 u + C_2 \theta)\,.
\end{equation} 
With $p=p(\theta)$ satisfying \eqref{adjointeq}, and $u=u(\theta)$ (also inserted into \eqref{adjointeq}) satisfying the state equation \eqref{modutheta0}, the derivative with respect to $\theta$ gives us the gradient of the reduced cost function $j_\alpha(\theta)=J_\alpha(u(\theta),\theta))$, which results in
\begin{align*}
 j_\alpha'(\theta)=& \frac{d}{d\theta} J_\alpha(u(\theta),\theta)) 
=\frac{d}{d\theta} \mathcal{L}_\alpha(u(\theta),\theta,p(\theta))) \\
&=\frac{\partial \mathcal{L}_\alpha}{\partial \theta} (u(\theta), \theta, p(\theta))
= \alpha I_R^X\left( C_2^* u + C_3 \theta \right) + L^\star p(\theta)\,.
\end{align*}
Now we can write down a new Lagrange functional for the computation of the Hessian matrix using the concept of second order adjoints. This new functional denoted by $\mathcal{L}^{\mathcal{H}}$ is given by
\begin{align*}
& \mathcal{L}^{\mathcal{H}}(u,\theta ,p; \hat{u}, \hat{\theta}, \hat{p}) :=  
\frac{\partial \mathcal{L}_\alpha}{\partial p} (u, \theta, p)\hat{p}
+\frac{\partial \mathcal{L}_\alpha}{\partial u} (u, \theta, p)\hat{u}
+\frac{\partial \mathcal{L}_\alpha}{\partial \theta} (u, \theta, p)\hat{\theta}\\
&=\mod(u,\theta)\hat{p}
+\langle \obs u-y,\obs\hat{u}\rangle_Y + \alpha \langle C_1 u + C_2 \theta,\hat{u}\rangle_U  
- \langle K\hat{u},p\rangle_{W',W}\\
&\qquad+ \alpha \langle C_2^* u + C_3 \theta,\hat{\theta}\rangle_X  
+ \langle L\hat{\theta},p\rangle_{W',W}.
\end{align*}
Again, the derivatives are computed and finally, the Hessian of the reduced cost functional can be given.
\begin{align*}
 j_\alpha''(\theta)(\hat{\theta},h)=& 
\frac{d^2}{d\theta^2} J_\alpha(u(\theta),\theta))(\hat{\theta},h)
=\frac{d}{d\theta} \mathcal{L}^H(u(\theta),\theta,p(\theta); \hat{u}, \hat{\theta}, \hat{p}) h\\
=&\frac{\partial \mathcal{L}^H}{\partial u} (u(\theta),\theta,p(\theta);\hat{u},\hat{\theta},\hat{p}) \frac{\partial u}{\partial\theta}(\theta) h
+\frac{\partial \mathcal{L}^H}{\partial \theta} (u(\theta),\theta,p(\theta);\hat{u},\hat{\theta},\hat{p}) h\\
&+\frac{\partial \mathcal{L}^H}{\partial p} (u(\theta),\theta,p(\theta);\hat{u},\hat{\theta},\hat{p}) \frac{\partial p}{\partial\theta}(\theta) h
\end{align*}
Choosing $\hat{u}=\hat{u}(\theta,\hat{\theta})$, $\hat{p}=\hat{p}(\theta,\hat{\theta})$ such that the first and the last term vanish for all $h$, i.e.,
\begin{align*}
 0 &=  \frac{\partial \mathcal{L}^H}{\partial u} (u(\theta),\theta,p(\theta);\hat{u},\hat{\theta},\hat{p}) \\
&=I_R^U \left(-K^*\hat{p}(\theta,\hat{\theta})+\obs^*\obs\hat{u}(\theta,\hat{\theta})+\alpha (C_1^*\hat{u}(\theta,\hat{\theta})+C_2\hat{\theta})\right),\\
 0 &=\frac{\partial \mathcal{L}^H}{\partial p} (u(\theta),\theta,p(\theta);\hat{u},\hat{\theta},\hat{p})
= -K\hat{u}(\theta,\hat{\theta})+L\hat{\theta}=0,
\end{align*}
where the lower equation is solved first and the resulting $\hat{u}(\theta)$ is inserted into the upper equation, which is the resolved for $\hat{p}(\theta)$.
This together with 
\[
\frac{\partial \mathcal{L}^H}{\partial \theta} (u,\theta,p;\hat{u},\hat{\theta},\hat{p}) h 
= \langle L h,\hat{p}\rangle_{W',W}+\alpha\langle C_2 h,\hat{u}\rangle_U+\alpha\langle C_3 h,\hat{\theta}\rangle_U
\]
yields, for any $h\in X$,
\[
 j_\alpha''(\theta)(\hat{\theta},h)
=\langle L^*\hat{p}(\theta,\hat{\theta})+\alpha\left(C_2^*\hat{u}(\theta,\hat{\theta})+C_3\hat{\theta}\right), h\rangle_X,
\]
thus, abbreviating $K^{-*}=(K^*)^{-1}$, 
\[
l j_\alpha''(\theta)(\hat{\theta},\check{\theta})
= \langle \left(L^*K^{-*} [\obs^*\obs +\alpha C_1]K^{-1}L+\alpha C_3
+\alpha L^*K^{-*}C_2+\alpha C_2^*K^{-1}L\right)\hat{\theta},\check{\theta}\rangle_X.
\]
Thus, altogether the gradient and Hessian of the reduced cost function $j$ read as follows
\begin{eqnarray*}
 j_\alpha'(\theta)
=\alpha I_R^X\left( C_2^* K^{-1}L + C_3 \right)\theta + L^\star K^{-*}\bigl[\obs^*(\obs K^{-1}L\theta - y) + \alpha (C_1 K^{-1}L + C_2) \theta\bigr]\,,\\
 j_\alpha''(\theta)=I_R^XH_\theta=I_R^X\Bigl(L^*K^{-*} [\obs^*\obs +\alpha C_1]K^{-1}L+\alpha C_3
+\alpha L^*K^{-*}C_2+\alpha C_2^*K^{-1}L\Bigr).
\end{eqnarray*}

If $L$ is invertible, then the role of $u$ and $\theta$ in the above procedure can be exchanged and for the reduced cost function $\ell_\alpha(u)=J_\alpha(u,\theta(u) )$ with $\theta(u)=L^{-1}(Ku-f)$
and $p(u)=-\alpha L^{-*}(C_2^*u+C_3\theta(u))$ we get 
\begin{align*}
& \ell_\alpha'(u)= \frac{d}{du} J_\alpha(u,\theta(u))) 
=\frac{d}{du} \mathcal{L}_\alpha(u,\theta(u),p(u)) 
=\frac{\partial \mathcal{L}_\alpha}{\partial u} (u, \theta(u), p(u))\\
 &= I_R^U\Bigl(\obs^*(\obs u-y)+\alpha \Bigl(C_1 u+ C_2 L^{-1}(Ku-f)+ K^* L^{-*}(C_2^* u + C_3 L^{-1}(Ku-f))\Bigr)\Bigr)
\end{align*}
and, with $\hat{\theta}(u,\hat{u})= L^{-1}K\hat{u}$, $\hat{p}(u,\hat{u})=-\alpha L^{-*}(C_2^*+C_3L^{-1}K)\hat{u}$,
\begin{align*}
 \ell_\alpha''(u)(\hat{u},\check{u})
&= \frac{d}{du} \mathcal{L}^H(u,\theta(u),p(u);\hat{u},\hat{\theta}(u,\hat{u}),\hat{p}(u,\hat{u}))\check{u}\\ 
&= \frac{\partial \mathcal{L}^H}{\partial u}(u,\theta(u),p(u);\hat{u},\hat{\theta}(u,\hat{u}),\hat{p}(u,\hat{u}))\check{u}\\
&=\langle -K^*\hat{p}(u,\hat{u})+\obs^*\obs\hat{u}+\alpha(C_1\hat{u}+C_2\hat{\theta}(u,\hat{u})),\check{u}\rangle_U\\
&=\langle \left[\obs^*\obs+\alpha(C_1+C_2L^{-1}K+K^*L^{-*}C_2^*+K^*L^{-*}C_3L^{-1}K)\right]\hat{u}
,\check{u}\rangle_U
\end{align*}
that is, $\ell_\alpha$ is the quadratic functional $\ell_\alpha(u)=\frac12\langle Hu,u\rangle +\langle g,u\rangle$ with
\begin{eqnarray*}
 &&
H = (I_R^U)^{-1}\ell_\alpha''(u) = \obs^*\obs+\alpha(C_1+C_2L^{-1}K+K^*L^{-*}C_2^*+K^*L^{-*}C_3L^{-1}K)\,,\\
 &&
g = (I_R^U)^{-1}\ell_\alpha'(0)= - \obs^*y-\alpha(C_2+K^*L^{-*}C_3)L^{-1}f.
\end{eqnarray*}

The computed Hessians can now be used to evaluate the MAP estimator for the posterior mean and the posterior covariance.

\subsubsection{Inverse source problem.}
$U=H_0^1(\Omega) \cap H^2(\Omega)$, $X=Y=W=W'=L^2(\Omega)$,\\ 
$\mod(u,\theta)=-\calA u+\theta$, $\obs u= u$, 
$K=\calA:H_0^1(\Omega) \cap H^2(\Omega)\to L^2(\Omega)$, $K^*=\calA^{-1}:L^2(\Omega)\to H_0^1(\Omega) \cap H^2(\Omega)$, 
$L=\mbox{id}:L^2(\Omega)\to L^2(\Omega)$, $L^*=\mbox{id}:L^2(\Omega)\to L^2(\Omega)$, 
$\obs:H_0^1(\Omega) \cap H^2(\Omega)\to L^2(\Omega)$, $\obs^*=\calA^{-2}:L^2(\Omega)\to H_0^1(\Omega) \cap H^2(\Omega)$, 
which results in
\begin{align*}
 g_{\theta}(\theta) &:=  (I_R^X)^{-1} j_\alpha'(\theta) \\
 &= \alpha  (C_2^*\calA^{-1} + C_3) \theta
+ \calA\bigl[\calA^{-2}(\calA^{-1}\theta-y) +\alpha(C_1\calA^{-1}+C_2)\theta\bigr],
\\
 H_{\theta}  &:=  (I_R^X)^{-1} j_\alpha'' = \calA(\calA^{-2} + \alpha C_1) \calA^{-1} + \alpha C_3 + \alpha \calA C_2 + \alpha C_2^* \calA^{-1},
\end{align*}
for the parameter $\theta$, where the computation of the adjoint is as described in Section~\ref{sec:adj}. 
For the state $u$ we have
\begin{align*}
 g_u(u) &:= (I_R^U)^{-1} \ell_\alpha'(u) \\
 &= \calA^{-2}(u-y) +\alpha \Bigl(C_1 u+ C_2 (\calA u-f)+ \calA^{-1} (C_2^* u + C_3 (\calA u-f)),\\
  H_u &:= (I_R^U)^{-1} \ell_\alpha''= \calA^{-2} + \alpha (C_1 + C_2 \calA + \calA^{-1} C_2^* + \calA^{-1}C_3 \calA)).
\end{align*}

\subsubsection{Backwards heat problem.}
$U=U_0 = \{ w \in L^2(H_0^1(\Omega)) \cap H^1(H^{-1}(\Omega)): w(0,x) = 0 \}$,
$W=L^2(H_0^1(\Omega))$,
$W'=L^2(H^{-1}(\Omega))$,
$X=H_0^1(\Omega)$,
$Y=L^2(\Omega)$,\\
$\mod(u,\theta)=(\partial_t+\calA) u + \ulI\calA\theta$, $\obs = \delta_T$, 
$K=(\partial_t+\calA):U_0\to W'$, $K^{-1}=\int_0^. e^{-\calA(.-s)} \cdot(s) \mbox{ d}s:W'\to U_0$,
$K^*=\int_0^. e^{-\calA(.-s)} \cdot(s) \mbox{ d}s:W'\to U_0$, $K^{-*}=K$,
$L=\ulI\calA:H_0^1(\Omega)\to W'$, 
$L^*=\calA^{-1}\int_0^T\cdot(s) \mbox{ d}s:W'\to H_0^1(\Omega)$, 
$K^{-1}L=\ulI-e^{\calA \cdot}:H_0^1(\Omega)\to U_0$, 
$L^*K^{-*}=\calA^{-1}\delta_T-\ulI^*:U_0\to H_0^1(\Omega)$, 
$\obs=\delta_T:U_0\to L^2(\Omega)$, 
$\obs^*= \frac{1}{2} (e^{-\calA(T-.)} - e^{-\calA (T+.)}): L^2(\Omega)\to U_0$, 
resulting in the following instances
\begin{align*}
& g_{\theta}(\theta) :=  (I_R^X)^{-1} j_\alpha'(\theta) = 
\alpha ( C_2^* (\ulI-e^{\calA \cdot}) + C_3)\theta + (\calA^{-1}\delta_T-\ulI^*) \\
 &\qquad \cdot \Bigl[\tfrac{1}{2} (e^{-\calA(T-.)} - e^{-\calA (T+.)})(\delta_T[\ulI-e^{\calA \cdot}]\theta - y) + \alpha (C_1 [\ulI-e^{\calA \cdot}] + C_2) \theta\Bigr],\\
& H_{\theta}  :=  (I_R^X)^{-1} j_\alpha'' = 
(\calA^{-1}\delta_T-\ulI^*)\tfrac{1}{2} (e^{-\calA(T-.)} - e^{-\calA (T+.)})(\delta_T[\ulI-e^{\calA \cdot}]\\
 &\qquad+\alpha \Bigl((\calA^{-1}\delta_T-\ulI^*)C_1 [\ulI-e^{\calA \cdot}] + (\calA^{-1}\delta_T-\ulI^*)C_2 + C_2^* (\ulI-e^{\calA \cdot}) + C_3  ,
\Bigr)
\end{align*}
for the parameter. 

Since $L$ is not invertible here, we substitute it by its right inverse $L^\sharp=\frac{1}{T}L^*$ to obtain, for the state,
\begin{align*}
& g_u(u) := (I_R^U)^{-1} \ell_\alpha'(u)=
\tfrac{1}{2} (e^{-\calA(T-.)} - e^{-\calA (T+.)}) ( \delta_T u -y) \\
 &\qquad + \alpha\Bigl(C_1u + C_2 \calA^{-1} \tfrac{1}{T}\tint_0^T((\partial_t+\calA)u-f)(s) \mbox{ d}s \\
 &\qquad+ \tfrac{1}{T} (\mbox{id}-e^{-\calA t})
\Bigl(C_2^*u + C_3 \tfrac{1}{T}\tint_0^T((\partial_t+\calA)u-f)(s) \mbox{ d}s\Bigr)\Bigr),
\\
&  H_u := (I_R^U)^{-1} \ell_\alpha'' =
\tfrac{1}{2}(e^{-\calA (T-.)} - e^{-\calA (T+.)}) \delta_T \\
 &\qquad + \alpha\Bigl(C_1 + C_2 \calA^{-1} \tfrac{1}{T}\tint_0^T(\partial_t+\calA)\cdot(s) \mbox{ d}s \\
 &\qquad+ \tfrac{1}{T} (\mbox{id}-e^{-\calA t})
\Bigl(C_2^* + C_3 \tfrac{1}{T}\tint_0^T(\partial_t+\calA)\cdot(s) \mbox{ d}s\Bigr).
\end{align*}

\subsection{Implementation.}
As already mentioned in the beginning of this section the analysis is done in python, with the help of Fenics and hippylib. With hippylib, the machinery of Bayesian inversion can be applied almost automatically to all stationary problems. However, concerning time-dependent problems a bit more work is necessary to get the results which are needed. In our case this means that the backwards heat problem had to be formulated in hippylib fashion to hand over to the algortihms for reconstruction. These algorithms were adapted to the all-at-once setting. In the following two sections the implementation of the problems as well as the results of the reconstruction are described in detail. 

From now on, we fix the spatial domain $\Omega$ on which both inverse problems are considered to be the unit square and space discretization is done by means of Lagrange finite elements. For time discretization in the backwards heat problem we use zero order discontinuos Galerkin method as already described previously, which leads to a backward Euler time stepping scheme. The discretized versions of paramter and state will be denoted by $u_h$ $\theta_h$, the corresponding finite dimensional function spaces by $U_h$ and $X_h$. 

To generate the synthetic observations the forward model is solved with the true parameter. Here, the dimension of the spatial discretization space for the state and the parameter is both 1681. In order to avoid an inverse crime, this forward simulation is done on a finer grid as compared to the reconstruction, where the dimension is chosen to be 961.
Discrete observation points are generated as 100 random points on the computational mesh, and the data at these points is perturbed with random noise. The measurement noise is assumed to follow a Gaussian distribution with zero mean and covariance $\Sigma = \delta^2 I$, where $\delta = 0.01$, and $I \in \mathbb{R}^{100\times100}$.
These synthetic observations are used (together with the prior) as input for the reconstructions.

\subsubsection{Inverse source problem.}
\begin{figure}
\centering
\begin{subfigure}[c]{0.3\textwidth}
\includegraphics[width=0.9\textwidth]{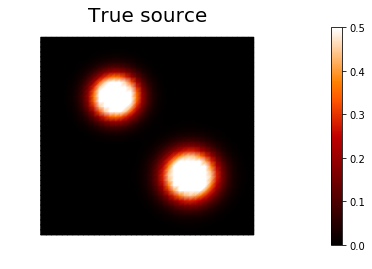}
\subcaption{True source.}
\label{is:truem}
\end{subfigure}
\begin{subfigure}[c]{0.3\textwidth}
\includegraphics[width=0.9\textwidth]{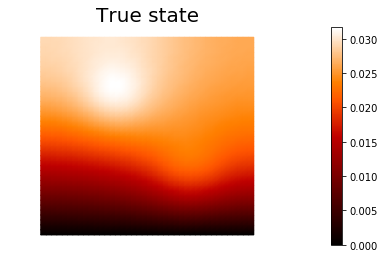}
\subcaption{True state.}
\label{is:trueu}
\end{subfigure}
\begin{subfigure}[c]{0.3\textwidth}
\includegraphics[width = 0.9\textwidth]{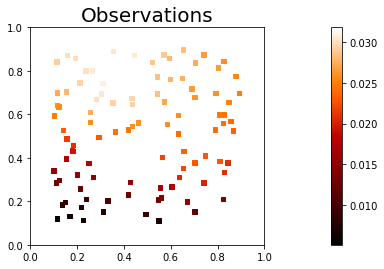}
\subcaption{Discrete observations.}
\label{is:obs}
\end{subfigure}
\caption{Simulated observations for the inverse source problem.}
\end{figure}

\begin{figure}
\centering
\begin{subfigure}[c]{0.35\textwidth}
\includegraphics[width=0.9\textwidth]{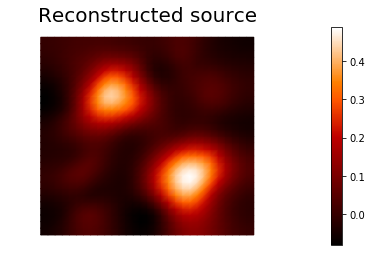}
\subcaption{Reconstructed source.}
\label{is:rm}
\end{subfigure}
\begin{subfigure}[c]{0.35\textwidth}
\includegraphics[width=0.9\textwidth]{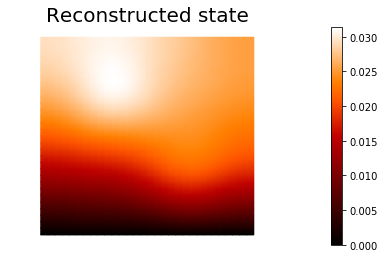}
\subcaption{Reconstructed state.}
\label{is:ru}
\end{subfigure}
\caption{Reconstructions for the inverse source problem.\label{is:rmu}}
\end{figure}

\begin{figure}
\centering
\begin{subfigure}[c]{0.47\textwidth}
\includegraphics[width=0.95\textwidth]{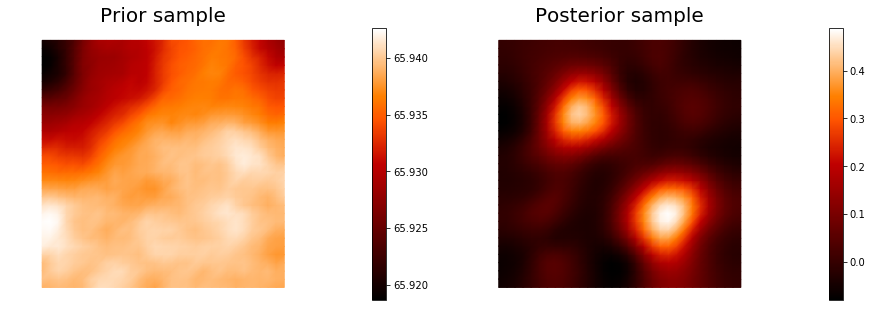}
\end{subfigure}
\begin{subfigure}[c]{0.47\textwidth}
\includegraphics[width=0.95\textwidth]{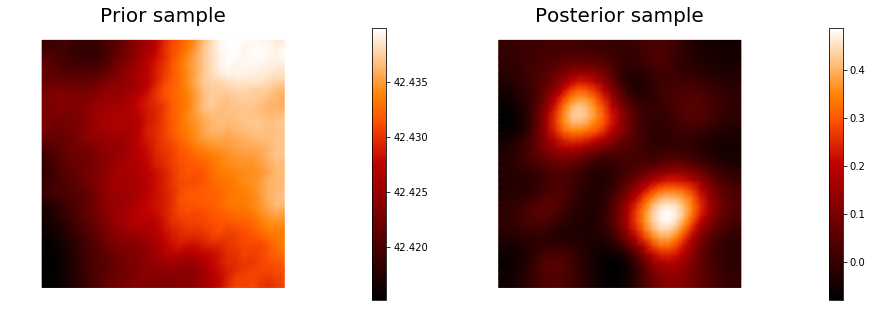}
\end{subfigure}
\caption{Samples from the posterior for the parameter variable of the inverse source problem.}
\label{is:samplespostpar}
\end{figure}

\begin{figure}
\centering
\begin{subfigure}[c]{0.47\textwidth}
\includegraphics[width=0.95\textwidth]{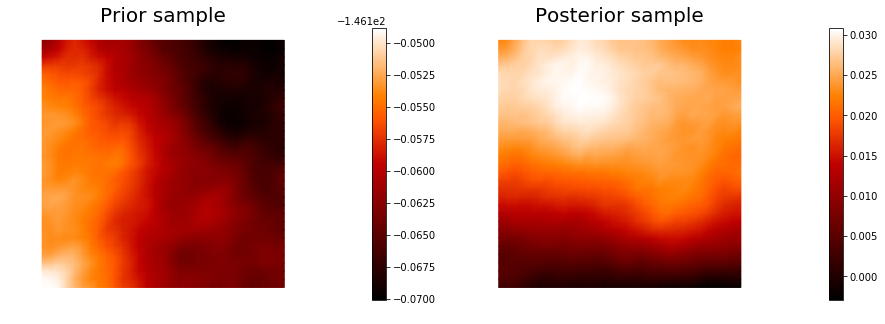}
\end{subfigure}
\begin{subfigure}[c]{0.47\textwidth}
\includegraphics[width=0.95\textwidth]{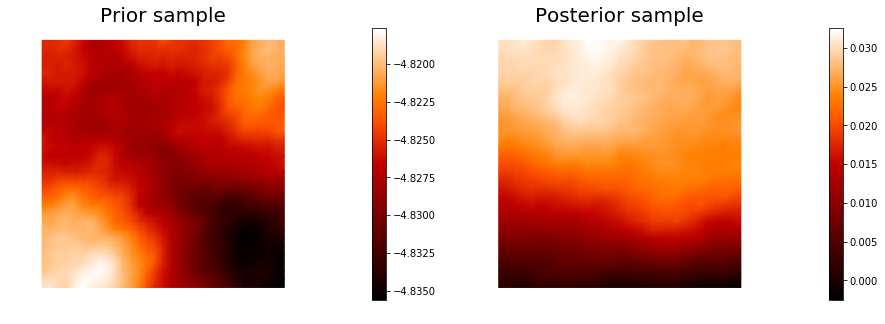}
\end{subfigure}
\caption{Samples from the posterior for the state variable of the inverse source problem.}
\label{is:samplespostst}
\end{figure}

\begin{figure}
\centering
\begin{subfigure}[c]{0.23\textwidth}
\includegraphics[width=0.95\textwidth]{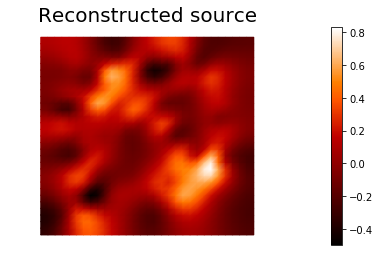}
\subcaption{Source: $\gamma_p = 1$ and $\kappa_p=10^{-2}$.}
\end{subfigure}
\begin{subfigure}[c]{0.23\textwidth}
\includegraphics[width=0.95\textwidth]{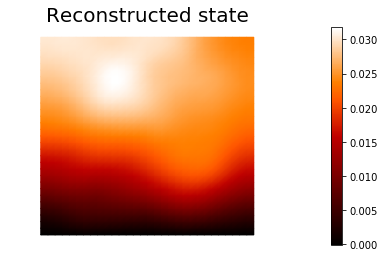}
\subcaption{State: $\gamma_p = 1$ and $\kappa_p=10^{-2}$.}
\end{subfigure}
\begin{subfigure}[c]{0.23\textwidth}
\includegraphics[width=0.95\textwidth]{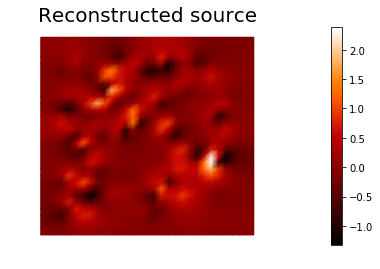}
\subcaption{Source: $\gamma_p = 10^{-2}$ and $\kappa_p=70$.}
\end{subfigure}
\begin{subfigure}[c]{0.23\textwidth}
\includegraphics[width=0.95\textwidth]{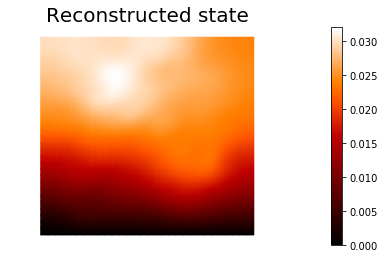}
\subcaption{State: $\gamma_p = 10^{-2}$ and $\kappa_p=70$.}
\end{subfigure}
\caption{Reconstructions for the inverse source problem for different values for $\kappa_p$ and $\gamma_p$ and $\kappa_s = 10^{-2}$ and $\gamma_s = 35$.}
\label{is:comparison}
\end{figure}

\begin{figure}
\centering
\begin{subfigure}[c]{0.23\textwidth}
\includegraphics[width=0.95\textwidth]{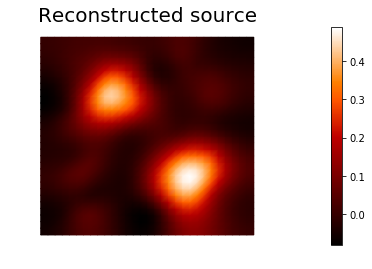}
\subcaption{Source: $\gamma_s = 10^{-2}$ and $\kappa_s=70$.}
\end{subfigure}
\begin{subfigure}[c]{0.23\textwidth}
\includegraphics[width=0.95\textwidth]{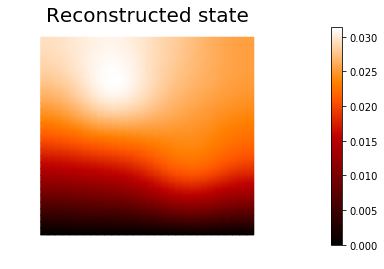}
\subcaption{State: $\gamma_s = 10^{-2}$ and $\kappa_s=70$.}
\end{subfigure}
\begin{subfigure}[c]{0.23\textwidth}
\includegraphics[width=0.95\textwidth]{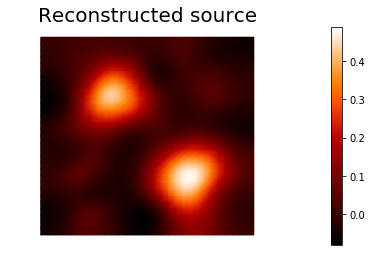}
\subcaption{Source: $\gamma_s = 10^{-5}$ and $\kappa_s=100$.}
\end{subfigure}
\begin{subfigure}[c]{0.23\textwidth}
\includegraphics[width=0.95\textwidth]{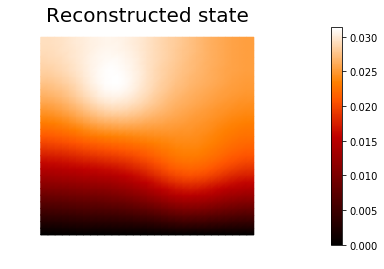}
\subcaption{State: $\gamma_s = 10^{-5}$ and $\kappa_s=100$.}
\end{subfigure}
\caption{Reconstructions for the inverse source problem for different values for $\kappa_s$ and $\gamma_s$ and $\kappa_p = 10^{-2}$ and $\gamma_p = 35$.}
\label{is:comparison2}
\end{figure}

The first test example considered for the implementation of the all-at-once version of Bayesian reconstruction is the inverse source problem. For this experiment the source is given by the function visualized in Figure \ref{is:truem}. As described above, the Poisson equation with this source is solved and random discrete observations are constructed from the resulting state, see Figure \ref{is:trueu} and \ref{is:obs}. To reconstruct the source these observations are used together with the prior distribution. In this case the prior distribution is defined by 
\begin{equation}
\label{prioris}
 \mu_{\mbox{prior}} = \mathcal{N} ( \mathbf{m}, \mathcal{C}_0), \mbox{ with } \mathbf{m} = \left( \begin{array}{c} m_u \\ m_{\theta} \end{array} \right),  \mathcal{C}_0 = \left( \begin{array}{cc} \kappa_s M + \gamma_s K & 0 \\ 0 & \kappa_p M + \gamma_p K \end{array} \right)^{-1},
\end{equation}
with $M$ and $K$ denoting the mass and the stiffness matrix, respectively, cf. \eqref{M},  \eqref{K}, and chosen parameters $m_{\theta}=0, \kappa_p = \kappa_s = 10^{-2}$ and $\gamma_p = \gamma_s = 35$. The prior matrix is chosen to be block diagonal, as operators on the off diagonal have not improved the reconstruction for this problem. The mean $m_{\theta}$ is chosen to be zero and the mean $m_u$ is the output of the forward problem when the source is $m_{\theta}$. The reconstruction using the synthetic observations and prior with the help of the CG-solver for the Hessian shown in the above computations leads to the results displayed in Figure \ref{is:rm} for the parameter and Figure \ref{is:ru} for the state. The posterior distribution is visualized with the help of samples for both the parameter in Figure \ref{is:samplespostpar} and the state in Figure \ref{is:samplespostst}.
Let us comment on the choice of the parameters $\kappa$ and $\gamma$ in the prior covariance. The reconstruction of the parameter appears to be strongly influenced by the choice of $\kappa_p$, $\gamma_p$ for the parameter prior, while the choice of $\kappa_s$, $\gamma_s$ in the state prior do not really influence the reconstruction neither for the parameter nor for the state. Some comparison of different choices of parameters can be seen in Figure~\ref{is:comparison} for the parameter and Figure~\ref{is:comparison2}. 
Nonzero choices for the covariance operators $C_2$ or $C_2^*$ other than zero have not improved the reconstructions for this example.

\begin{figure}
\centering
\begin{subfigure}[c]{0.47\textwidth}
\includegraphics[width=0.95\textwidth]{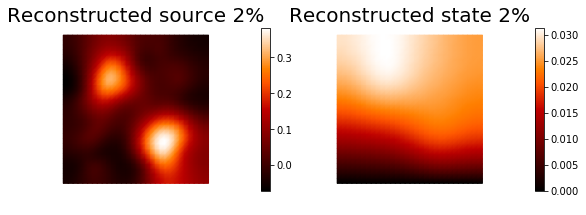}
\end{subfigure}
\begin{subfigure}[c]{0.47\textwidth}
\includegraphics[width=0.95\textwidth]{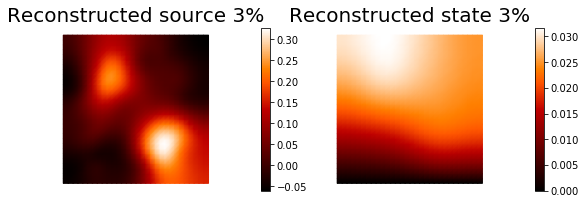}
\end{subfigure}
\caption{Reconstructions for the inverse source problem with different noise levels.}
\label{is:noiselevel1}
\end{figure}

\begin{figure}
\centering
\begin{subfigure}[c]{0.47\textwidth}
\includegraphics[width=0.95\textwidth]{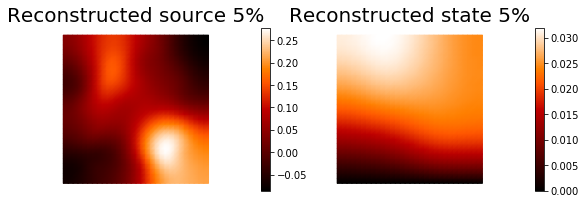}
\end{subfigure}
\begin{subfigure}[c]{0.47\textwidth}
\includegraphics[width=0.95\textwidth]{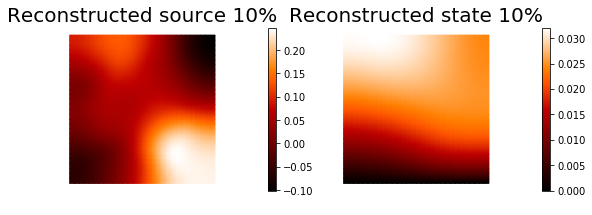}
\end{subfigure}
\caption{Reconstructions for the inverse source problem with different noise levels.}
\label{is:noiselevel2}
\end{figure}

In Figures \ref{is:rmu} -- \ref{is:comparison2}, the noise was set to 1\% as described in the introduction of this section. However, the reconstructions are quite good up to a noise level of 3\%, as it can be seen in Figure~\ref{is:noiselevel1} and Figure~\ref{is:noiselevel2}.

\subsubsection{Backwards heat equation, sampled initial condition.}\label{sec:bh_sampinit}
For the backwards heat equation the prior as described in Section \ref{sec:priors} is used, with the covariance operator for the initial condition chosen to be
\begin{equation}
C_p = (\kappa M + \gamma K)^{-1},
\end{equation}
where $M$ and $K$ are the mass and stiffness matrices resulting from the finite element discretization, with $\kappa=1.5$ and $\gamma=0.5$. The time variable is discretized in the same way as in Section \ref{sec:eig} resulting in the covariance operator
\begin{equation}
C_s = \mbox{diag}((e^{-t_i K} C_p e^{-t_i K})_{i=0}^N)
\end{equation} 
for the state variable,
so that altogether we arrive at 
\[
\mathcal{C}_0 = \left( \begin{array}{cc}C_s & 0 \\ 0 & C_p \end{array} \right).
\] 
Also here, nonzero covariance operators $C_2$ or $C_2^*$ did not improve the reconstructions but would only complicate computation of powers of $\mathcal{C}_0$.

\begin{figure}
\centering
\begin{subfigure}[c]{0.3\textwidth}
\includegraphics[width=0.9\textwidth]{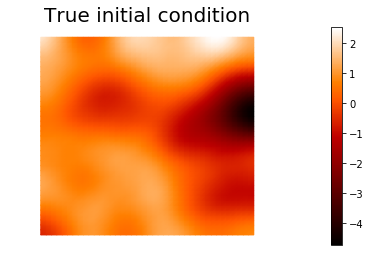}
\subcaption{True initial condition.}
\label{bh:truem}
\end{subfigure}
\begin{subfigure}[c]{0.3\textwidth}
\includegraphics[width=0.9\textwidth]{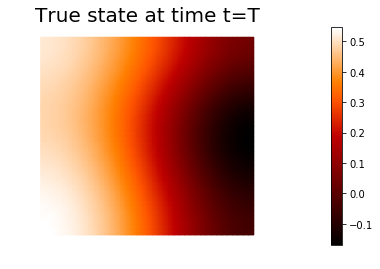}
\subcaption{True state at time $T$.}
\label{bh:trueu}
\end{subfigure}
\begin{subfigure}[c]{0.3\textwidth}
\includegraphics[width = 0.9\textwidth]{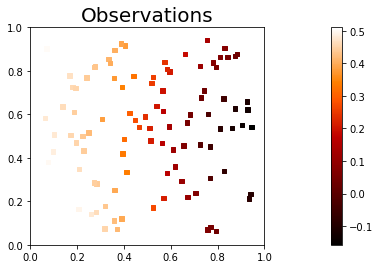}
\subcaption{Discrete observations.}
\label{bh:obs}
\end{subfigure}
\caption{Simulated observations for the backwards heat problem.}
\end{figure}

\begin{figure}
\centering
\begin{subfigure}[c]{0.35\textwidth}
\includegraphics[width=0.9\textwidth]{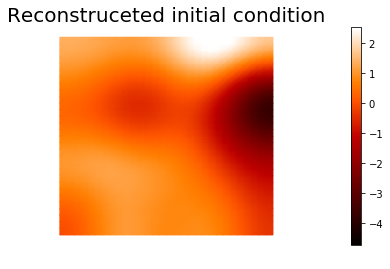}
\end{subfigure}
\begin{subfigure}[c]{0.35\textwidth}
\includegraphics[width=0.9\textwidth]{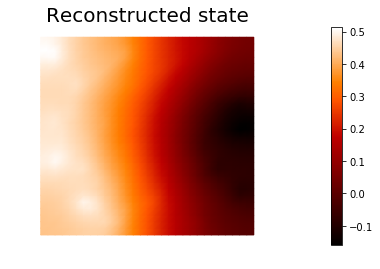}
\end{subfigure}
\caption{Reconstructions for the backwards heat problem.}
\label{bh:reconstr}
\end{figure}

\begin{figure}
\centering
\begin{subfigure}[c]{0.47\textwidth}
\includegraphics[width=0.95\textwidth]{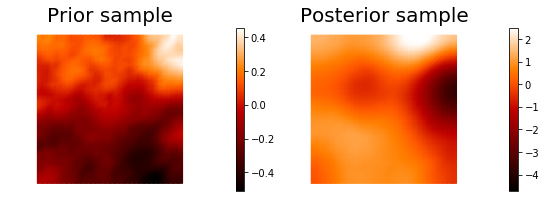}
\end{subfigure}
\begin{subfigure}[c]{0.47\textwidth}
\includegraphics[width=0.95\textwidth]{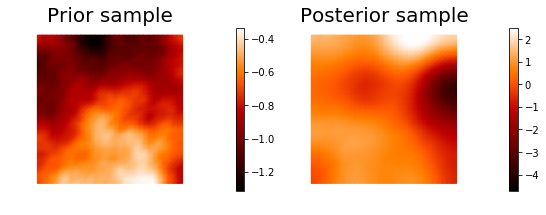}
\end{subfigure}
\caption{Samples from the posterior for the initial condition of the backwards heat problem.}
\label{bh:samplespostpar}
\end{figure}

Here, we have the setting $T=0.1, N=4$ and $t_i = i \frac{T}{N}$. The prior mean for the inital condition is chosen to be zero, therefore, also the mean for the state is zero for each $t \in [0,T]$.
Due to the diagonal structure of the operator matrix it is not necessary to save the whole matrix, instead for every time step only the application of the matrix vector product is needed. This overcomes the problem of block matrix usage in python. For the experiment, the initial condition is sampled from the prior distribution as seen in  Figure \ref{bh:trueu}. Then, the forward problem is solved (Figure \ref{bh:trueu} and discrete synthetic observations of the state at time $T=0.1$ are constructed, see Figure \ref{bh:obs}. Finally, the MAP estimators can be computed as described in the beginning of this section for both the parameter and the state variable, see Figure \ref{bh:reconstr}. In this example the posterior samples are only shown for the parameter in Figure \ref{bh:samplespostpar}.

\subsubsection{Backwards heat equation, chosen initial condition.}
\begin{figure}
\centering
\begin{subfigure}[c]{0.3\textwidth}
\includegraphics[width=0.9\textwidth]{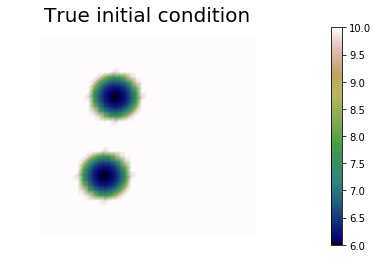}
\subcaption{True initial condition.}
\label{bh:truem2}
\end{subfigure}
\begin{subfigure}[c]{0.3\textwidth}
\includegraphics[width=0.9\textwidth]{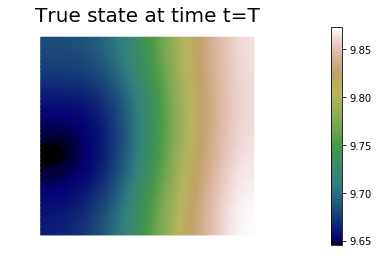}
\subcaption{True state at time $T$.}
\label{bh:trueu2}
\end{subfigure}
\begin{subfigure}[c]{0.3\textwidth}
\includegraphics[width = 0.9\textwidth]{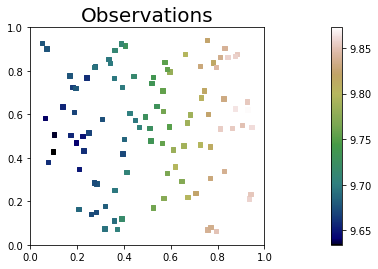}
\subcaption{Discrete observations.}
\label{bh:obs2}
\end{subfigure}
\caption{Simulated observations for the backwards heat problem.}
\end{figure}

\begin{figure}
\centering
\begin{subfigure}[c]{0.35\textwidth}
\includegraphics[width=0.9\textwidth]{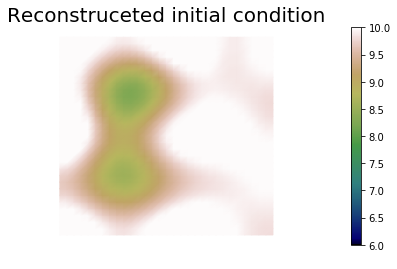}
\end{subfigure}
\begin{subfigure}[c]{0.35\textwidth}
\includegraphics[width=0.9\textwidth]{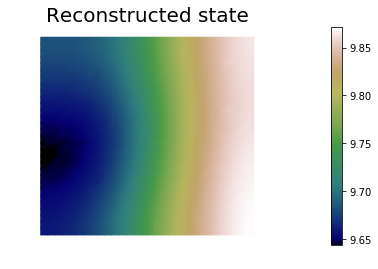}
\end{subfigure}
\caption{Reconstructions for the backwards heat problem.}
\label{bh:r2}
\end{figure}

\begin{figure}
\centering
\begin{subfigure}[c]{0.47\textwidth}
\includegraphics[width=0.95\textwidth]{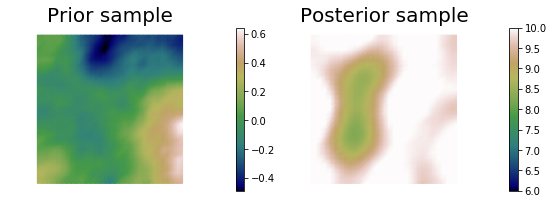}
\end{subfigure}
\begin{subfigure}[c]{0.47\textwidth}
\includegraphics[width=0.95\textwidth]{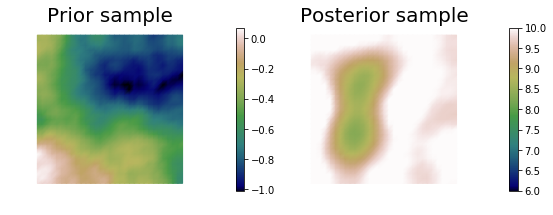}
\end{subfigure}
\caption{Samples from the posterior for the initial condition of the backwards heat problem.}
\label{bh:samplespostpar2}
\end{figure}

We now demonstrate performance of the method with a fixed initial condition as visualized in Figure~\ref{bh:truem2}. Again, discrete synthetic observations are made which are used for the reconstruction, see Figure~\ref{bh:obs2}. The prior is chosen as in Subsection \ref{sec:bh_sampinit}. 
The reconstructions of the parameter and the state are shown in Figure~\ref{bh:r2}. Samples are shown in Figure \ref{bh:samplespostpar2}.

\subsubsection{Backwards heat equation, chosen initial condition, prior motivated by the link condition.} \label{sec:heuristic}
\begin{figure}
\centering
\begin{subfigure}[c]{0.35\textwidth}
\includegraphics[width=0.9\textwidth]{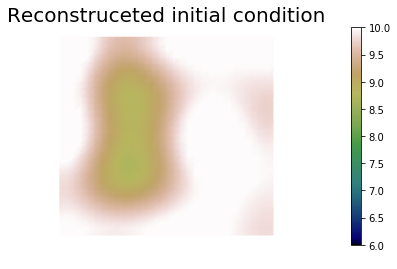}
\end{subfigure}
\begin{subfigure}[c]{0.35\textwidth}
\includegraphics[width=0.9\textwidth]{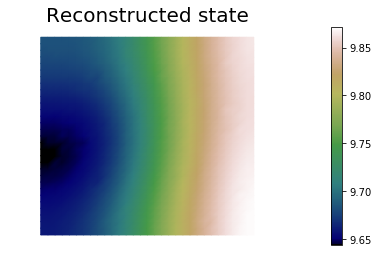}
\end{subfigure}
\caption{Reconstructions for the backwards heat problem with the heuristic prior.}
\label{bh:heuristic}
\end{figure}

Referring to Section \ref{sec:heuristicprior} we do the same experiments as in the last section with the prior given by $\tilde{\mathcal{C}}$ as in \eqref{Ctil} but skipping the off-diagonal blocks for ease of implementation
\begin{equation}
 \mu_{\mbox{prior}} = \mathcal{N} ( \mathbf{m}, \mathcal{C}_0), \mbox{ with } \mathbf{m} = \left( \begin{array}{c} m_u \\ m_{\theta} \end{array} \right), \quad  \mathcal{C}_0 = \left( \begin{array}{cc} \mathcal{M} & 0 \\ 0 & K^{1/2} \end{array} \right).
\end{equation}
The operator $\mathcal{M}$ denotes an operator consisting of the mass matrix $M$ for every timestep $t \in [0,T]$. The setup of the numerical experiment is the same as in the last example concerning discretization, initial condition and observations. The reconstructions can be seen in Figure~\ref{bh:heuristic}, which suggests quite a good fit also for this prior.

\section{Conclusions and Remarks.}
In this paper we have combined the Bayesian approach with an all-at-once formulation of the inverse problem. We have done so for linear problems, in particular focusing on two prototypical examples, namely the inverse source problem for the Poisson equation and the backwards heat equation.
\\
Our next step will be to extend this approach to nonlinear problems such as the identification of coefficients in time dependent and stationary PDEs. Moreover we will further investigate the use of joint state and parameter priors, possibly also taking covariances between them into account.

\section*{Acknowledgment} This work was supported by the Austrian Science
Fund FWF under the grants P30054 and DOC 78.

\end{document}